\RequirePackage{etoolbox}
\csdef{input@path}{{style/}{graphics/}}
\documentclass[aos,MSNbibl,seceqn,dvips]{arximspdf}
\usepackage{mathrsfs}
\usepackage{graphicx}

\doi{10.1214/12-AOS1036} 
\volume{40}
\issue{4}
\pubyear{2012}
\firstpage{2266}
\lastpage{2292}

\makeatletter
\newcommand{\rrvert}{\vert}
\newcommand{\llvert}{\vert}
\newcommand{\eqref}[1]{(\ref{#1})}
\newcommand{\V}[1]{\mathbf{#1}}
\newcommand{\VV}[1]{\bolds{#1}}
\newtheorem{thmm}{Theorem}[section]
\newtheorem{lemma}{Lemma}[section]
\newtheorem{prop}{Proposition}[section]
\newtheorem{cor}{Corollary}[section]
\newcommand{\argmax}{\operatorname{arg\,max}}
\makeatother

\begin{document}
\begin{frontmatter}

\title{Consistency of community detection in networks under
degree-corrected stochastic block models}
\runtitle{Consistency of community detection}

\begin{aug}
\author[1]{\fnms{Yunpeng} \snm{Zhao}\ead[label=e1]{yzhao15@gmu.edu}},
\author[2]{\fnms{Elizaveta} \snm{Levina}\corref{}\ead[label=e2]{elevina@umich.edu}\thanksref{t1}}
\and
\author[2]{\fnms{Ji} \snm{Zhu}\ead[label=e3]{jizhu@umich.edu}\thanksref{t2}}
\thankstext{t1}{Supported in part by NSF Grants DMS-08-05798, DMS-01106772
and DMS-1159005}
\thankstext{t2}{Supported in part by NSF Grant
DMS-07-48389 and NIH Grant R01-GM-096194.}

\runauthor{Y. Zhao, E. Levina and J. Zhu}
\affiliation{George Mason University, University of Michigan and
University of Michigan}
\address[1]{Y. Zhao\\
Department of Statistics\\
George Mason University\\
4400 University Drive, MS 4A7\\
Fairfax, Virginia 22030-4444\\
USA\\
\printead{e1}}

\address[2]{E. Levina\\
J. Zhu\\
Department of Statistics\\
University of Michigan\\
439 West Hall \\
1085 S. University Ave. \\
Ann Arbor, Michigan 48109-1107 \\
USA\\
\printead{e2}\\
\phantom{E-mail:\ }\printead*{e3}}
\end{aug}

\received{\smonth{11} \syear{2011}}
\revised{\smonth{7} \syear{2012}}

%
\begin{abstract}
Community detection is a fundamental problem in network analysis, with
applications in many diverse areas. The stochastic block model is a
common tool for model-based community detection, and asymptotic tools
for checking consistency of community detection under the block model
have been recently developed.
However, the block
model is limited by its assumption that all nodes within a community
are stochastically equivalent, and provides a poor fit to networks with
hubs or highly varying node degrees within communities, which are
common in practice. The degree-corrected stochastic block model
was proposed to address this shortcoming and allows
variation in node degrees within a community while preserving the
overall block community structure. In this paper we establish general
theory for checking consistency of community detection under the
degree-corrected stochastic block model and compare several community
detection criteria under both the standard and the degree-corrected
models. We show which criteria are consistent under which models and
constraints, as well as compare their relative performance in practice.
We find that methods based on the degree-corrected block model, which
includes the standard block model as a special case, are consistent
under a wider class of models and that modularity-type methods require
parameter constraints for consistency, whereas likelihood-based methods
do not. On the other hand, in practice, the degree correction involves
estimating many more parameters, and empirically we find it is only
worth doing if the node degrees within communities are indeed highly
variable. We illustrate the methods on simulated networks and on a
network of political blogs.
\end{abstract}

%

\end{frontmatter}
%
\section{Introduction}
Networks have become one of the more common forms of data, and network
analysis has received a lot of attention in computer science, physics,
social sciences, biology and statistics (see \cite{Getoor2005,Goldenberg2010,Newman2010} for reviews). The applications are many and
varied, including social networks \cite{Wasserman1994,robinsetc2007}, gene regulatory networks \cite{Schlitt2007},
recommender systems and security monitoring. One of the fundamental
problems in network analysis is community detection, where communities
are groups of nodes that are, in some sense, more similar to each other
than to other nodes. The precise definition of community, like that of
a cluster in multivariate analysis, is difficult to formalize, but many
methods have been developed to address this problem (see \cite
{Goldenberg2010,NewmanPNAS,Fortunato2010} for comprehensive recent
reviews), often relying on the intuitive notion of community as a group
of nodes with many links between themselves and fewer links to the rest
of the network.

Three groups of methods for community detection can be loosely
identified in the literature. A number of greedy algorithms such as
hierarchical clustering have been proposed (see \cite{Newman2004Review}
for a review), which we will not focus on in this paper. The second
class of methods involves optimization of some ``reasonable'' global
criteria over all possible network partitions and includes graph cuts
[\cite*{Shi00}, \cite*{Wei&Cheng1989}], spectral clustering \cite{Ng01} and
modularity \cite{Newman&Girvan2004,NewmanPNAS}, the latter discussed
in detail below. Finally, model-based methods rely on fitting a
probabilistic model for a network with communities. Perhaps the best
known such model is the stochastic block model, which we will also
refer to as simply the block model \mbox{\cite
{Holland83,Snijders&Nowicki1997,Nowicki2001}}. Other models include a
recently introduced degree-corrected stochastic block model \cite
{Karrer10}, mixture models for directed networks \cite
{Newman&Leicht2007}, multivariate latent variable models \cite
{Handcock2007}, latent feature models \cite{Hoff2007} and mixed
membership stochastic block models for modeling overlapping communities
\cite{Airoldi2008}. From the algorithmic point of view, many
model-based methods also lead to criteria to be optimized over all
partitions, such as the profile likelihood under the assumed model.

The large number of available methods leads to the question of how to
compare them in a principled manner, other than on individual examples.
There has been little theoretical analysis of community detection
methods until very recently, when a consistency framework for community
detection was introduced by Bickel and Chen \cite{Bickel&Chen2009}.
They developed general theory for checking the consistency of detection
criteria under the stochastic block model (discussed in detail below)
as the number of nodes grows and the number of communities remains
fixed, and their result has been generalized to allow the number of
communities to grow in \cite{Choietal2011}; see also \cite{Rohe2011}.
The stochastic block model, however, has serious limitations in
practice: it treats all nodes within a community as stochastically
equivalent, and thus does not allow for the existence of ``hubs,''
high-degree nodes at the center of many communities observed in real
data. To address this issue, Karrer and Newman \cite{Karrer10} proposed
the degree-corrected stochastic block model, which can accommodate hubs
(a similar model for a directed network was previously proposed in
\cite
{wang1987}, but they did not focus on community detection and assumed
known community membership). In~\cite{Karrer10}, the authors gave
several examples showing this model fits data with hubs much better
than the block model; however, there are no consistency results
available under this new model, and thus no way to compare methods in general.

In this paper we generalize the consistency framework of \cite
{Bickel&Chen2009} to the degree-corrected stochastic block model and
obtain a general theorem for community detection consistency. Since the
degree-corrected model includes the regular block model as a special
case, consistency results under the block model follow automatically.
We then evaluate two types of modularity and the two criteria derived
from the block model and the degree-corrected block model using this
general framework. One of our goals is to emphasize the difference
between assumed models (needed for theoretical analysis)
and criteria for finding the optimal partition, which may or may not be
motivated by a particular model. What we ultimately show agrees with
statistical common sense: criteria derived from a particular model are
consistent when this model is assumed, but not necessarily consistent
if the model does not hold. Further, if a criterion relies implicitly
on an assumption about the model parameters (e.g., modularity
implicitly assumes that links within communities are stronger than
between), then it will be consistent only if the model parameters are
constrained to satisfy this assumption. We make all of the above
statements precise later in the paper.

The rest of the article is organized as follows. We set up all notation
and define the relevant models and criteria in Section \ref
{secmodels}. Consistency results under the regular and the
degree-corrected stochastic block models for all of the criteria in
Section \ref{secmodels} are stated in Section \ref{seccriteria}. The
general consistency theorem which implies all of these results is
presented in Section \ref{secthm}. In Section \ref{secsim} we compare
the performance of these criteria on simulated networks, and in Section
\ref{secdata} we illustrate the methods on a network of political
blogs. Section \ref{secsummary} concludes with a summary and
discussion. All proofs are given in the \hyperref[app]{Appendix}.

\section{Network models and community detection criteria}
\label{secmodels}
Before we proceed to discuss specific criteria and models, we introduce
some basic notation. A network $N=(V,E)$, where $V$ is the set of nodes
(vertices), $|V| = n$, and $E$ is the set of edges, can be represented
by its $n\times n$ adjacency matrix $A=[A_{ij}]$, where $A_{ij}=1$ if
there is an edge from $i$ to $j$, and $A_{ij}=0$ otherwise. We only
consider unweighted and undirected networks here, and thus $A$ is a
binary symmetric matrix. The community detection problem can be
formulated as finding a disjoint partition $V = V_1 \cup\cdots\cup
V_K$ or, equivalently, a set of node labels $\V{e}=\{e_1,\ldots,e_n\}$,
where $e_i$ is the label of node $i$ and takes values in $\{1,2,\ldots,K\}$.

For any set of label assignments $\V{e}$, let $O(\V{e})$ be the
$K\times K$ matrix defined by
\[
O_{kl}(\V{e})=\sum_{ij} A_{ij}I
\{e_i=k,e_j=l\} ,
\]
where $I$ is the indicator function. Further, let
\[
O_{k}(\V{e})=\sum_{l} O_{kl}(
\V{e}) ,\qquad  L=\sum_{ij} A_{ij} .
\]
For $k \neq l$, $O_{kl}$ is the total number of edges between
communities $k$ and $l$; $O_{k}$ is the sum of node degrees in
community $k$, and $L$ is the sum of all degrees in the network.
If self-loops are not allowed (i.e., $A_{ii} = 0$ is enforced), then we
can also interpret $O_{kk}$ as twice the total number of edges within
community $k$ and $L$ as twice the number of edges in the whole
network. Finally, let $n_k(\V{e})=\sum_i I\{e_i=k\}$ be the number of
nodes in the $k$th community, and $f(\V{e})= (\frac{n_1}{n},
\frac
{n_2}{n},\ldots,\frac{n_K}{n}  )^T$.

The stochastic block model, which is perhaps the most commonly used
model for networks with communities, postulates that, given node labels
$\V{c} = \{c_1, \ldots, c_n\}$, the edge variables $A_{ij}$'s are
independent Bernoulli random variables with
%
\begin{equation}
\label{bm} E[A_{ij}] = P_{c_i c_j} ,
\end{equation}
where $P=[P_{ab}]$ is a $K \times K$ symmetric matrix. We will use this
formulation throughout the paper, which allows for self-loops. While it
is also common to exclude self-loops, sometimes they are present in the
data (as in our example in Section \ref{secdata}) and allowing them
leads to simpler notation. In principle, all of our results go through
for the version of the models with self-loops excluded, with
appropriate modifications made to the proofs.

Under the model \eqref{bm}, all nodes with the same label are
stochastically equivalent to each other, which in practice limits the
applicability of the stochastic block model, as pointed out
in \cite{Karrer10}. The alternative proposed in \cite{Karrer10}, the
degree-corrected stochastic block model, is to replace \eqref{bm} with
%
\begin{equation}
\label{dcbm} E[A_{ij}] = \theta_i \theta_j
P_{c_i c_j} ,
\end{equation}
where $\theta_i$ is a ``degree parameter'' associated with node $i$,
reflecting its individual propensity to form ties. The degree
parameters have to satisfy a constraint to be identifiable, which in
\cite{Karrer10} was set to $\sum_i \theta_i I(c_i = k) = 1$, for each
$k$ (other constraints are possible). Further, they replaced the
Bernoulli likelihood by the Poisson, to simplify technical derivations.
With these assumptions, a profile likelihood can be derived by
maximizing over $\theta$ and $P$, giving the following criterion to be
optimized over all possible partitions:
%
\begin{equation}
\label{dcbm-lik} Q_\mathrm{DCBM}(\V{e}) = \sum_{kl}
O_{kl}\log\frac{O_{kl}}{O_k O_l} .
\end{equation}
We have compared the performance of this criterion in practice to its
slightly more complicated version based on the (correct) Bernoulli
likelihood instead of the Poisson and found no difference in the
solutions these two methods produce. The Bernoulli distribution with a
small mean is well approximated by the Poisson distribution, and most
real networks are sparse, so one can expect the approximation to work
well; see also a more detailed discussion of this in \cite{Perry2012}.
We will use~\eqref{dcbm-lik} in all further analysis, to be consistent
with \cite{Karrer10} and take advantage of the simpler form.

The degree-corrected model includes the regular stochastic block model
as a special case, with all $\theta_i$'s equal. Enforcing this
additional constraint on the profile likelihood leads to the following
criterion to be optimized over all partitions:
%
\begin{equation}
\label{bm-lik} Q_\mathrm{BM}(\V{e}) = \sum_{kl}
O_{kl}\log\frac{O_{kl}}{n_k n_l} .
\end{equation}
Like criterion \eqref{dcbm-lik}, this is based on the Poisson
assumption but gives identical results to the Bernoulli version in
practice. Here we use the form \eqref{bm-lik} for consistency with
\eqref{dcbm-lik} and with \cite{Karrer10}.

A different type of criterion used for community detection is
modularity, introduced in \cite{Newman&Girvan2004}; see also \cite
{NewmanPNAS} and \cite{Newman2006}. The basic idea of modularity is to
compare the number of observed edges within a community to the number
of expected edges under a null model and maximize this difference over
all possible community partitions. Thus, the general form of a
modularity criterion is
%
\begin{equation}
\label{generalM} Q(\V{e})=\sum_{ij}[A_{ij}-P_{ij}]I(e_i
= e_j),
\end{equation}
where $P_{ij}$ is the (estimated) probability of an edge falling
between $i$ and $j$ under the null model. The convention in the physics
literature is to divide $Q$ by $L$, which we omit here, since it does
not change the solution.

The choice of the null model, that is, of a model with no communities
($K=1$), determines the exact form of modularity. The stochastic block
model with $K=1$ is simply the Erdos--Renyi random graph, where $P_{ij}$
is a constant which can be estimated by $L/n^2$. Plugging $P_{ij} =
L/n^2$ into \eqref{generalM} gives what we will call the Erdos--Renyi
modularity (ERM),
%
\begin{equation}
\label{erm} Q_\mathrm{ERM}(\V{e})=\sum_{k}
\biggl(O_{kk}-\frac{n_k^2}{n^2}L \biggr) .
\end{equation}

\begin{table}[b]
\caption{Summary of community detection criteria}
\label{tablecrit}
\begin{tabular*}{\textwidth}{@{\extracolsep{\fill}}lcc@{}}
\hline
& \textbf{Block model} & \multicolumn{1}{c@{}}{\textbf{Degree-corrected block model}} \\
\hline
Modularity & $\sum_{k} (O_{kk}-\frac{n_k^2}{n^2}L)$ (ERM) & $ \sum_{k} ( O_{kk}-\frac{O_{k}^2}{L^2} L) $ (NGM)
\\[3pt]
Likelihood & $\sum_{kl} O_{kl}\log\frac{O_{kl}}{n_kn_l}$ (BM) & $\sum_{kl} O_{kl}\log\frac{O_{kl}}{O_k O_l}$ (DCBM) \\
\hline
\end{tabular*}
\end{table}

If instead we take the degree-corrected model with $K=1$ as the null
model, it postulates that $P_{ij} \propto\theta_i \theta_j$, where
$\theta_i$ is the degree parameter. This is essentially the well-known
expected degree random graph, also known as the configuration model. In
this case, $P_{ij}$ can be estimated by $d_i d_j/L$, where $d_i = \sum_j A_{ij}$
is the degree of node $i$. Substituting this into \eqref
{generalM} gives the popular Newman--Girvan modularity (NGM), introduced
in \cite{Newman&Girvan2004}:
%
\begin{equation}
\label{ngm} Q_\mathrm{NGM}(\V{e})=\sum_{k}
\biggl( O_{kk}-\frac{O_{k}^2}{L^2} L\biggr).
\end{equation}

The four different criteria for community detection are summarized in
Table~\ref{tablecrit}. Note that the two likelihood-based criteria, BM
and DCBM, take into account all links within and between communities,
and which communities they connect; whereas the modularities would not
change if all the links connecting different communities were randomly
permuted (as long as they did not become links within communities).
Further, note that the degree correction amounts to substituting $O_k$
for $n_k$ and $L$ for $n$, both for modularity and likelihood-based
criteria. Thus, if all nodes within a community are treated as
equivalent, their number suffices to weigh community strength
appropriately; and if the nodes are allowed to have different expected
degrees, then the number of edges becomes the correct weight. Both of
these features make sense intuitively and, as we will see later, will
fit in naturally with consistency conditions.

Our analysis indicates that Newman--Girvan modularity and
degree-corrected block model criteria are consistent under the more
general degree-corrected models but Erdos--Renyi modularity and block
model criteria are not, even though they are consistent under the
regular block model. Further, we show that likelihood-based methods are
consistent under their assumed model with no restrictions on
parameters, whereas modularities are only consistent if the model
parameters are constrained to satisfy a ``stronger links within than
between'' condition, which is the basis of modularity derivations. In
short, we show that a criterion is consistent when the underlying model
and assumptions are correct, and not necessarily otherwise.

\section{Consistency of community detection criteria}
\label{seccriteria}

Here we present all the consistency results for the four different
criteria defined in Section \ref{secmodels}. All these results follow
from the general consistency theorem in Section \ref{secthm}; the
proofs are given in the \hyperref[app]{Appendix}. The notion of consistency of
community detection as the number of nodes grows was introduced in
\cite
{Bickel&Chen2009}. They defined a community detection criterion $Q$ to
be consistent if the node labels obtained by maximizing the criterion,
$\hat{\V{c}} = \argmax_{\V{e}} Q(\V{e})$, satisfy
%
\begin{equation}
P[\hat{\V{c}}=\V{c}] \rightarrow1 \qquad \mbox{as } n \rightarrow \infty .
\label{strong}
\end{equation}
Strictly speaking, this definition suffers from an identifiability
problem, since most reasonable criteria, including all the ones
discussed above, are invariant under a permutation of community labels
$\{1, \ldots, K\}$. Thus, a better way to define consistency is to
replace the equality $\hat c = c$ with the requirement that $\hat c$
and $c$ belong to the same equivalence class of label permutations. For
simplicity of notation, we still write $\hat c = c$ in all consistency
results in the rest of the paper, but take them to mean that $\hat c$
and $c$ are equal up to a permutation of labels.

The notion of consistency in \eqref{strong} is very strong, since it
requires asymptotically no errors. One can also define what we will call
weak consistency,
%
\begin{equation}
\label{weak} \forall\varepsilon>0\qquad  \mathbb{P} \Biggl[ \Biggl(\frac{1}{n}
\sum_{i=1}^n 1(\hat{c}_i \neq
c_i) \Biggr) < \varepsilon \Biggr] \rightarrow 1 \qquad\mbox{as }n
\rightarrow\infty,
\end{equation}
where equality is also interpreted to mean membership in the same
equivalence class with respect to label permutations. In \cite
{Bickel&Chen2012}, conditions were established for a criterion to be
weakly consistent under the stochastic block model. All other
assumptions being equal, weak consistency only requires that the
expected degree of the graph $\lambda_n \rightarrow\infty$, whereas
strong consistency requires $\lambda_n / \log n \rightarrow\infty$.
Here, we will analyze both strong and weak consistency under the
degree-corrected stochastic block model.

For the asymptotic analysis, we use a slightly different formulation
of the degree-corrected model than that given by \cite{Karrer10}. The
main difference is that we treat true community labels $\V{c}$ and
degree parameters $\VV{\theta} = (\theta_1, \ldots, \theta_n)$ as latent
random variables rather than fixed parameters. Note, however, that the
criteria we analyze were obtained as profile likelihoods with
parameters treated as constants. This is one of the standard approaches
to random effects models, known as conditional likelihood (see page~234 of
\cite{McCulloch01}). The network model we use for consistency analysis
can be described as follows:
\begin{longlist}[(1)]
\item[(1)] Each node is independently assigned a pair of latent variables
$(c_i, \theta_i)$, where $c_i$ is the community label taking values in
$1, \ldots, K$,
and $\theta_i$ is a discrete ``degree variable'' taking values in $x_1
\leq\cdots\leq x_M$. We do not assume that $c_i$ is independent of
$\theta_i$.
\item[(2)] The marginal distribution of $c$ is multinomial with parameter
$\VV{\pi}=(\pi_1, \ldots, \pi_K)^T$, and $\theta$ satisfies $E[\theta_i]=1$
for identifiability.
\item[(3)] Given $\V{c}$ and $\VV{\theta}$, the edges $A_{ij}$ are
independent Bernoulli random variables with
\[
E[A_{ij}|\V{c},\VV{\theta}]=\theta_i\theta_j
P_{c_ic_j} ,
\]
where $P=[P_{ab}] $ is a $K \times K$ symmetric matrix.
\end{longlist}
For simplicity, we allow self-loops in the network, that is, $
E[A_{ii}|\V{c},\VV{\theta}]=\theta_i^2 P_{c_ic_i}$. Otherwise diagonal
terms of $A$ have to be treated separately, which ultimately makes no
difference for the analysis but makes notation more awkward.

To ensure that all probabilities are always less than 1, we require the
model to satisfy the constraint $x_M^2 \max_{a,b}P_{ab} \leq1$. We
also need to consider how the model changes with $n$. If $P_{ab}$
remains fixed as $n$ grows, the expected degree $\lambda_n$ will be
proportional to $n$, which makes the network unrealistically dense.
Instead, we allow the matrix $P$ to scale with $n$ and, in a slight
abuse of notation, reparameterize
it as $P_n=\rho_n P$, where $\rho_n = P(A_{ij} = 1) \rightarrow0$ and
$P$ is fixed. We then specify the rate of $\mathbf{c}$ the expected degree
$\lambda_n = n \rho_n$, which has to satisfy $\frac{\lambda_n}{\log n}
\rightarrow\infty$ for strong consistency and $\lambda_n \rightarrow
\infty$ for weak consistency.

Let $\Pi$ be the $K \times M$ matrix representing the joint
distribution of $(c_i,\theta_i)$ with $\mathbb{P}(c_i=a,\theta_i=x_u)=\Pi_{au}$. Further, define $\tilde{\pi}_a=\sum_u x_u \Pi_{au}$.
Note that $\sum_{a} \tilde{\pi}_a=1$ since $\mathbb{E}(\theta_i)=1$.
Moreover, we have $\tilde{\pi}_a=\pi_a$ if $\V{c}$ and $\VV{\theta
}$ are
independent, or if $\theta_i \equiv1$ (block models). Thus, we can
view $\tilde{\pi}$ as an adjusted version of $\pi$.


Next, we state our consistency results for the two types of
modularities under both the degree-corrected and the standard block model.

\begin{thmm} \label{corNGM}
Under the degree-corrected stochastic block model, if the parameters satisfy
\[
\tilde{\mathcal{E}}_{aa}>0 ,\qquad \tilde{\mathcal{E}}_{ab}<0\qquad
\mbox{for all } a\neq b,
\]
where $ \tilde{P}_0=\sum_{ab} \tilde{\pi}_a \tilde{\pi}_b P_{ab}$,
$\tilde{W}_{ab}=\frac{\tilde{\pi}_a \tilde{\pi}_b P_{ab}}{\tilde
{P}_0}$, $\tilde{\mathcal{E}}=\tilde{W}-(\tilde{W} \V{1})(\tilde
{W} \V{1})^T$,
the Newman--Girvan modularity is strongly consistent when $\lambda_n /
\log n \rightarrow\infty$ and weakly consistent when $\lambda_n
\rightarrow\infty$.
\end{thmm}
The parameter constraints in Theorem \ref{corNGM} require,
essentially, that the links within communities are more likely than the
links between. This is particularly easy to see when $K=2$, in which
case the constraint simplifies to
\[
P_{11}P_{22}>P_{12}^2.
\]

Taking $\theta_i \equiv1$, we immediately obtain the following.
%
\begin{cor}[(Established in \cite{Bickel&Chen2009})]
\label{corNGM-block} Under the standard stochastic block model with
parameters satisfying Theorem \ref{corNGM} constraints with $\tilde
\pi
$ replaced by~$\pi$, Newman--Girvan modularity is strongly consistent
when $\lambda_n / \log n \rightarrow\infty$ and weakly consistent when
$\lambda_n \rightarrow\infty$.
\end{cor}

For Erdos--Renyi modularity, which has not been studied theoretically
before, we can also show consistency under the standard block model,
albeit with a slightly stronger condition on links within communities
being more likely than the links between:
%
\begin{thmm}
\label{corERM}
Under the standard stochastic block model, if the parameters satisfy
\[
P_{aa}>P_0 ,\qquad P_{ab}<P_0 \qquad\mbox{for
all } a\neq b ,
\]
where $P_0=\sum_{ab} \pi_a \pi_b P_{ab}$,
the Erdos--Renyi modularity criterion \eqref{erm} is strongly consistent
when $\lambda_n / \log n \rightarrow\infty$ and weakly consistent when
$\lambda_n \rightarrow\infty$.
\end{thmm}

However, the Erdos--Renyi modularity is not consistent under the
degree-corrected model, at least not under the same parameter
constraint. The Erdos--Renyi modularity prefers to group nodes with
similar degrees together, which may not agree with true communities
when the variance in node degrees is large. Here is a counter-example
demonstrating this. Let $K=2, \VV{\pi}=(1/2,1/2)^T$, $\rho_n = 1$ (so
that the graph becomes dense as $n \rightarrow\infty$), and
\[
P=\pmatrix{ 0.1 & 0.05 \vspace*{2pt}
\cr
0.05 & 0.1 }.
\]
Further, $\VV{\theta}$ is independent of $\V{c}$ and takes only two
values, $1.6$ and $0.4$, with probability $1/2$ each. If we assign all
nodes their true labels, the population version of the criterion (where
all random quantities are replaced by their expectations under the true
model) gives $Q_\mathrm{ERM}= 0.0125$. However, by grouping nodes with the
same value of $\theta_i$'s together, we get the population version of
$Q_\mathrm{ERM} = 0.0135$, higher than the value for the true partition,
and this solution will therefore be preferred in the limit.


Once again, the result makes sense intuitively, since the Erdos--Renyi
modularity uses the regular block model as its null hypothesis, and the
parameter constraint matches the ``fewer links between than within''
notion. From the algorithmic point of view, the main difference between
Erdos--Renyi modularity and Newman--Girvan modularity is that the latter
depends on the edge matrix $O$ only and ``weighs'' communities by the
number of edges, whereas the former weighs communities by the number of
nodes $n_k$ (which, under the block model, is proportional to the
number of edges, but under the degree-corrected model is not).

Next we state the consistency results for the two criteria derived from
profile likelihoods, DCBM \eqref{dcbm-lik} and BM \eqref{bm-lik}. These
require no parameter constraints.

\begin{thmm}\label{corDCBL}
Under the degree-corrected stochastic block model (and therefore under
the regular model as well), the degree-corrected criterion~\eqref
{dcbm-lik} is strongly consistent when $\lambda_n / \log n \rightarrow
\infty$ and weakly consistent when $\lambda_n \rightarrow\infty$.
\end{thmm}

\begin{thmm}\label{corBL}
Under the stochastic block model, the block model criterion~\eqref
{bm-lik} is strongly consistent when $\lambda_n / \log n \rightarrow
\infty$ and weakly consistent when $\lambda_n \rightarrow\infty$.
\end{thmm}
Theorem \ref{corBL} was proved in \cite{Bickel&Chen2009} for a
slightly different form of the profile likelihood (Bernoulli rather
than the Poisson). Under the degree-corrected block model, criterion
\eqref{bm-lik} is not necessarily consistent---the same
counter-example can be used to demonstrate this. As was the case with
modularities, the criterion consistent under the degree-corrected block
model depends on $O$ only, whereas the criterion consistent only under
the regular block model also depends on $n_k$.

The theoretical results suggest that the likelihood-based criteria are
always preferable over the modularity-based criteria, and that criteria
based on the degree-corrected model are always preferred to the
criteria based on the regular block model, since they are consistent
under weaker conditions. In practice, however, this may not always
hold. Computationally, modularity type criteria can be approximately
optimized by solving an eigenvalue problem \cite{Newman2006}, whereas
likelihood type criteria have no such approximations and thus have to
be optimized by slower heuristic search algorithms, as was done in
\cite
{Bickel&Chen2009} and \cite{Karrer10}. Moreover, fitting the
degree-corrected block model requires estimating many more parameters
than fitting a block model and creates the usual trade-off between
model complexity and goodness of fit. If the node degrees within
communities do not vary widely, fitting a block model may provide a
better solution; see more on this in Section \ref{secsim}.

\section{A general theorem on consistency under degree-corrected
stochastic block models} \label{secthm}
Here we prove a general theorem for checking consistency under
degree-corrected stochastic block models for \textit{any} criterion
defined by a reasonably nice function. All consistency results for
specific methods discussed in Section \ref{seccriteria} are
corollaries of this theorem.

A large class of community detection criteria can be written as
%
\begin{equation}
\label{generalForm} Q(\V{e})=F \biggl( \frac{O(\V{e})}{\mu_n},f(\V{e}) \biggr) ,
\end{equation}
where $\mu_n=n^2 \rho_n$. For instance, many graph cut methods (mincut,
ratio cut \cite{Wei&Cheng1989}, normalized cut \cite{Shi00}) have this
form and use functions that are designed to minimize the number of
edges between communities. All criteria discussed in Section~\ref
{seccriteria} can also be written in
this form. Our goal here is to establish conditions for consistency of
a criterion of this form under degree-corrected block models.

A natural condition for consistency is that the ``population version''
of $Q(\V{e})$ should be maximized by the correct community assignment,
as in $M$-estimation. To define the population version of $Q$, we first
define functions $H(S)$ and $h(S)$ corresponding to population versions
of $O(\V{e})$ and $f(\V{e})$, respectively (the precise meaning of
``population version'' is clarified in Proposition \ref{prop1} below).
For any generic array $S=[S_{kau}] \in
\mathcal{R}^{K \times K \times M}$, define a $K \times K$ matrix
$H(S)=[H_{kl}(S)]$ by
\[
H_{kl}(S)=\sum_{abuv} x_u
x_v P_{ab} S_{kau} S_{lbv} ,
\]
and a $K$-dimensional vector $h(S)=[h_{k}(S)]$ by
\[
h_{k}(S)=\sum_{au} S_{kau} .
\]
Also define $R(\V{e}) \in\mathcal{R}^{K \times K \times M}$ by
\[
R_{kau}(\V{e})  =\frac{1}{n}\sum_{i=1}^n
I(e_i=k,c_i=a,\theta_i=x_u) .
\]
Then we have the following:
%
\begin{prop}{\label{prop1}}
%
\begin{eqnarray}
\frac{1}{\mu_n} \mathbb{E} [O_{kl}| \V{c},\VV{\theta} ] & =&
H_{kl}\bigl(R(\V {e})\bigr) , \label{prop1a}
\\
f_k(\V{e}) & =& h_k\bigl(R(\V{e})\bigr) .
\label{prop1b}
\end{eqnarray}
\end{prop}
%
Proposition \ref{prop1} explains the precise meaning of ``population
version'': we take the conditional expectations given $\V{c}$ and $\VV{\theta}$
and write them as functions of a generic variable $S$ instead
of $R(\V{e})$. The population version of $Q$ is defined as $F(H(S),h(S))$.

Now we can specify the key sufficient condition as follows:
\begin{itemize}[($\ast$)]
\item[($\ast$)] $F(H(S),h(S))$ is uniquely maximized over $\mathscr
{S}=\{S\dvtx  S \geq0, \sum_k S_{kau}=\Pi_{au} \}$ by $S=\mathbb{D}$, with
$\mathbb{D}_{kau}=\Pi_{au}E_{ka}$, for any $a$ and $u$, where $E$ is
any row permutation of a $K \times K$ identity matrix.
\end{itemize}
The matrix $E$ deals with the permutation equivalence class. Since
$R(\V
{c})\rightarrow\mathbb{D}$ as $n \rightarrow\infty$, $S=\mathbb{D}$
implies each class $k$ exactly matches a community in the population.
For simplicity, in what follows we assume that $E$ is in fact the
identity matrix itself. We will elaborate on this condition below. In
addition, we need some regularity conditions, analogous to those in
\cite{Bickel&Chen2009}:
\begin{longlist}[(a)]
\item[(a)] $F$ is Lipschitz in its arguments;
\item[(b)] Let $W=H(\mathbb{D})$. The directional derivatives $\frac
{\partial^2 F}{\partial\varepsilon^2}(M_0+\varepsilon(M_1-M_0), \V
{t}_0+\varepsilon(\V{t}_1-\V{t}_0))|_{\varepsilon=0+}$ are continuous in
$(M_1,\V{t}_1)$ for all $(M_0,\V{t}_0)$ in a neighborhood of $(W,\VV{\pi})$;
\item[(c)] Let $G(S)=F(H(S),h(S))$. Then on $\mathscr{S}$, $\frac
{\partial G((1-\varepsilon)\mathbb{D}+\varepsilon S)}{\partial\varepsilon
}|_{\varepsilon=0+}<-C<0$ for all $\VV{\pi},P$.
\end{longlist}

Now we are ready to state the main theorem.
%
\begin{thmm}\label{genThm}
For any $Q(\V{e})$ of the form \eqref{generalForm}, if $\VV{\pi},P,F$
satisfy ($\ast$), \textup{(a)--(c)}, then
$Q$ is strongly consistent under degree-corrected stochastic block
models if $\frac{\lambda_n}{\log n}\rightarrow\infty$ and weakly
consistent if $\lambda_n \rightarrow\infty$.
\end{thmm}
The proof is given in the \hyperref[app]{Appendix}. This theorem is a generalization of
Theorem~1 in \cite{Bickel&Chen2009} from the standard stochastic block
models to degree-corrected models, and it implies all of the
consistency results in Section~\ref{seccriteria}.

Finally, we return to the key condition ($\ast$).
If $Q(\V{e})$ is maximized by the true community labels $\V{c}$, then
as $n\rightarrow\infty$, $F(H(S),h(S))$, the population version of
$Q(\V{e})$, should also be maximized by the true partition $S=\mathbb
{D}$, since $R(\V{c})\rightarrow\mathbb{D}$ and $Q(\V{c})
\rightarrow
F(H(\mathbb{D}),h(\mathbb{D}))$, making ($\ast$) a natural condition.
Further, since for any $\V{e}$, $\sum_{k} R_{kau}(\V{e}) \rightarrow
\Pi_{au}$, the limit $S$ of $R(\V{e})$ must satisfy $\sum_k S_{kau}=\Pi_{au}$. Therefore, we only need to consider maximizers of
$F(H(S),h(S))$ satisfying this constraint.

\section{Numerical evaluation}
\label{secsim}

In this section we compare the performance of the four community
detection criteria from Section \ref{secmodels} on simulated data,
generated from the regular or the degree-corrected block model. The
criteria are maximized over partitions using a greedy label-switching
algorithm called tabu search \cite{Beasley1998,Glover&Laguna1997}. The
key idea of tabu search is that once a node label has been switched, it
will be ``tabu'' and not available for switching for a certain number
of iterations, to prevent being trapped in a local maximum. Even though
tabu search cannot guarantee convergence to the global maximum, it
performs well in practice. Moreover, we run the search for a number of
initial values and different orderings of nodes, to help avoid local maxima.

To compare the solution to the true labels, we use the adjusted Rand
index \cite{Hubert1985}, a measure of similarity between partitions
commonly used in clustering. We have also computed the normalized
mutual information, a measure more commonly used by physicists in the
networks literature, which gives very similar results (not reported to
save space). The adjusted Rand index is scaled so that 1 corresponds to
the perfect match and 0 to the expected difference between two random
partitions, with higher values indicating better agreement. The figures
in this section all present the median adjusted Rand index over 100
replications.

In all examples below, we generate networks with $n=1000$ nodes and
$K=2$ communities. The node labels are generated independently with
$P(c_i = 1) = \pi$, $P(c_i = 2) = 1 - \pi$. By varying $\pi$, we can
investigate robustness of the methods to unbalanced community sizes.
The probability matrix for the block model and the degree-corrected
block model is set to
\begin{eqnarray*}
P= \rho\pmatrix{ 4 & 1 \vspace*{2pt}
\cr
1 & 4 } ,
\end{eqnarray*}
where we vary $\rho$ to obtain different expected degrees $\lambda$.

\subsection{The degree-corrected stochastic block model} For this
simulation, we generate data from the degree-corrected model with two
possible values for the degree parameter $\theta$. The degree
parameters are generated independently from the labels, with
\[
P(\theta_{i}= m x) = P(\theta_{i}= x) = 1/2 ,
\]
which implies $x = \frac{2}{m+1}$, since we need to have $E(\theta_i) =
1$. We vary the ratio $m$ from~1 (the regular block model) to 10, which
allows us to study the effect of model misspecification on the regular
block model. In this simulation, the community sizes are balanced ($\pi
= 0.5$).

\begin{figure}[b]

\includegraphics{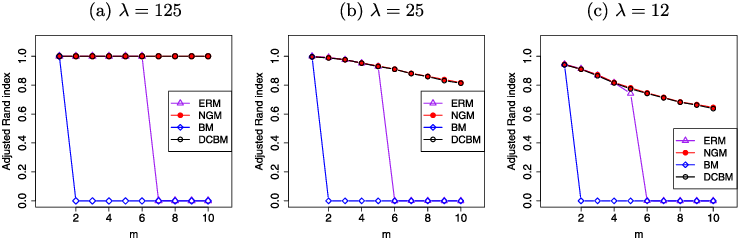}

\caption{Results for the degree-corrected stochastic block model with
two values for the degree parameters, $\pi=0.5$, $m$ varies.}
\label{figsim1}
\end{figure}

%

Figure \ref{figsim1} shows the results for three different expected
degrees $\lambda$. For the densest network with $\lambda= 125$ in
Figure \ref{figsim1}(a), the degree-corrected block model and
Newman--Girvan modularity perform the best overall, as they assume the
correct model and the methods are consistent. At $m=1$, the regular
block model is just as good, but its performance deteriorates rapidly
as $m$ increases. The Erdos--Renyi modularity also performs perfectly
for $m=1$, and it takes larger values of $m$ for its performance to
deteriorate than for block model likelihood, so we can conclude that
the Erdos--Renyi modularity is more robust to variation in degrees. For
both of them, poor results are due to grouping nodes with similar
degrees together. The overall trend for sparser networks [Figure \ref
{figsim1}(b) and (c)] is similar, but all methods perform worse, as
with fewer links there is effectively less data to use for fitting the
model, and the effect is more pronounced for large $m$, when degrees
have higher variance.

%
\begin{figure}

\includegraphics{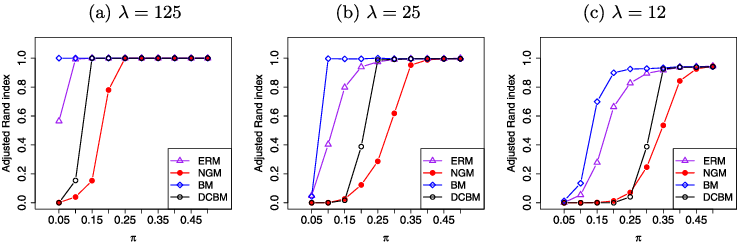}

\caption{Results for the standard stochastic block model, $m=1$, $\pi$ varies.}
\label{figsim2}
\end{figure}

%
\subsection{The stochastic block model} Here we focus on the standard
stochastic block model ($m=1$) and vary $\pi$ to assess robustness to
unbalanced community sizes. All the four criteria are consistent in
this case, but, in practice, the closer~$\pi$ is to 0.5, the better
they perform (Figure \ref{figsim2}), with the exception of the block
model likelihood in the dense case ($\lambda= 125$), where it performs
perfectly for all $\pi$. Overall, the block model likelihood performs
best, which is natural because it is the maximum likelihood estimator
of the correct model. The Erdos--Renyi modularity also performs better
than the other two criteria, which overfit the data by assuming the
degree-corrected model and accounting for variation in observed
degrees, which in this case only adds noise.
%

\subsection{Unbalanced community sizes}
%
\begin{figure}

\includegraphics{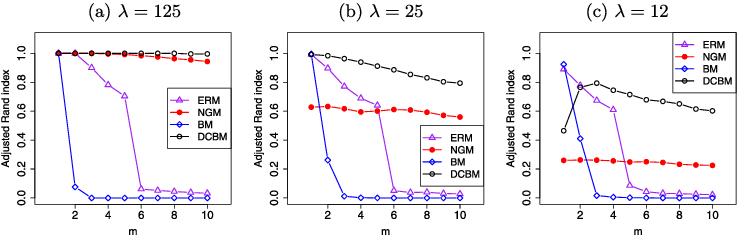}

\caption{Results for the degree-corrected stochastic block model with
two values for the degree parameters, $\pi=0.3$, $m$ varies.}
\label{figsim31}
\end{figure}
%
%
%

In this simulation we consider the degree-corrected stochastic block
model with unbalanced community sizes. We fix $\pi= 0.3$ and vary the
ratio $m$ in Figure \ref{figsim31}. For a dense network [$\lambda=
125$, Figure \ref{figsim31}(a)], the performance with $\pi=0.3$ is
similar to the balanced case with $\pi= 0.5$ [Figure~\ref
{figsim1}(a)]. However, in sparser networks modularity performs much
worse with unbalanced community sizes. This can also be seen in Figure
\ref{figsim2} for the case $m=1$. The failure of modularity to deal
with unbalanced community sizes was also recently pointed out by \cite
{Zhang&Zhao2012}. Note also that in the sparsest case ($\lambda= 12$,
Figure \ref{figsim31}), the degree-corrected model suffers from
over-fitting when $m=1$, as was also seen in Figure \ref{figsim2}.


\subsection{A different degree distribution}
In the last simulation we test the sensitivity of all methods, but in
particular the degree-corrected model, to the assumption of a discrete
degree distribution. Here we sample the degree parameters $\theta_i$
independently from the following distribution:
\[
\theta_i= \cases{ \eta_i,
& \quad $\mbox{w.p. $\alpha$,}$
\vspace*{2pt}\cr
2/(m+1), & \quad $\mbox{w.p. $(1-\alpha)/2$,}$
\vspace*{2pt}\cr
2m/(m+1), & \quad $\mbox{w.p. $(1-\alpha)/2$,}$ }
\]
where $\eta_i$ is uniformly distributed on the interval $[0,2]$. The
variance of $\theta_i$ is
equal to $\alpha/ 3 + (1-\alpha) (m-1)^2 / (m+1)^2$. In this
simulation, we fix $m = 10$, which makes the variance a decreasing
function of $\alpha$, and vary $\alpha$ from 0 to 1. We also fix $\pi=0.5$.

\begin{figure}[b]

\includegraphics{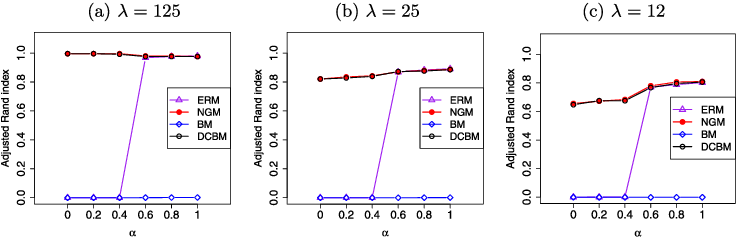}

\caption{Results for the degree-corrected stochastic block model with a
mixture degree distribution, $m = 10$, $\pi = 0.5$, mixture parameter $\alpha$ varies.}
\label{figsim4}
\end{figure}

%
The results in Figure \ref{figsim4} show that the degree-corrected
block model likelihood and Newman--Girvan modularity still perform well,
which suggests that the discreteness of $\theta$ is not a crucial
assumption. The regular block model fails in this case, as we would
expect from earlier results since $m=10$, but the performance of the
Erdos--Renyi modularity improves as $\alpha$ increases, which agrees
with our earlier observation on its relative robustness to variation in degrees.

\section{Example: The political blogs network}
\label{secdata}

In this section we analyze a real network of political blogs compiled
by \cite{Adamic05}. The nodes of this network are blogs about US
politics and the edges are hyperlinks between these blogs. The data
were collected right after the 2004 presidential election and
demonstrate strong divisions; each blog was manually labeled as liberal
or conservative by \cite{Adamic05}, which we take as ground truth.
Following the analysis in \cite{Karrer10}, we ignore directions of the
hyperlinks and focus on the largest connected component of this
network, which contains 1222 nodes, 16,714 edges and has the average
degree of approximately~27. Some summary statistics of the node degrees
are given in Table \ref{tabdegrees}, which shows that the degree
distribution is heavily skewed to the right.
%
\begin{table}
\caption{Statistics of node degrees in the political blogs network}
\label{tabdegrees}
\begin{tabular*}{\textwidth}{@{\extracolsep{\fill}}lccccc@{}}
\hline
\multicolumn{1}{@{}l}{\textbf{Mean}} & \textbf{Median} & \textbf{Min} & \textbf{1st Qt.} & \textbf{3rd Qt.} & \multicolumn{1}{c@{}}{\textbf{Max}} \\
\hline
27.36 & 13.00 & 1.00 & 3.00 & 36.00 & 351.00\\
\hline
\end{tabular*}
\end{table}


\begin{figure}

\includegraphics{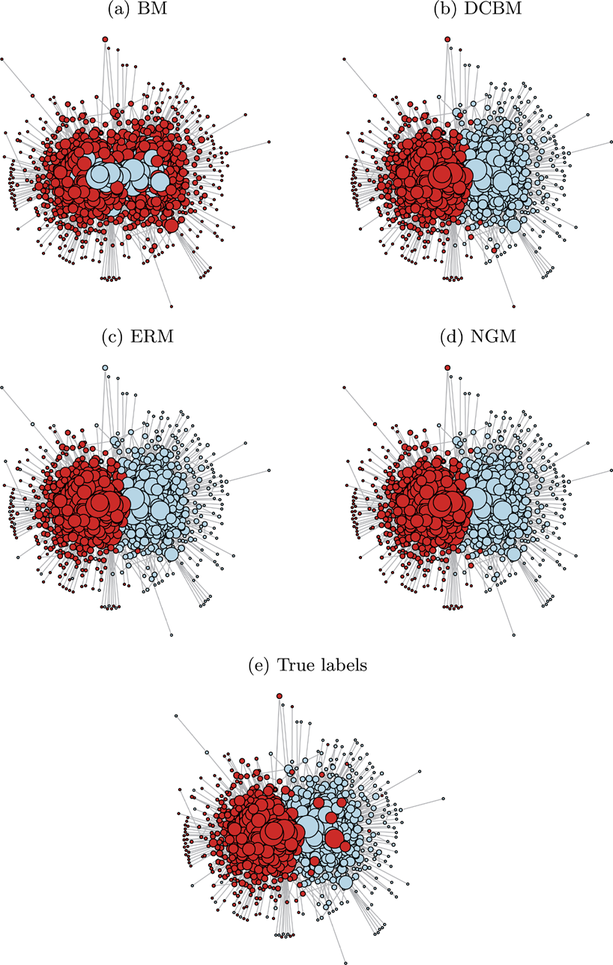}

\caption{Political blogs data. Node area is proportional to the
logarithm of its degree and the colors represent community labels.}
\label{figpolblog}
\end{figure}

We compare the partitions into two communities found by the four
different community detection criteria with the true labels using the
adjusted Rand index. The Newman--Girvan modularity and the
degree-corrected model find very similar partitions (they differ over
only four nodes and have the same adjusted Rand index value of 0.819,
the highest of all methods). The partition found by the Erdos--Renyi
modularity has a slightly worse agreement with the truth (adjusted Rand
index of 0.793). The block model likelihood divides the nodes into two
groups of low degree and high degree, with the adjusted Rand index of
nearly 0, which is equivalent to random guessing. The results are shown
in Figure \ref{figpolblog} (drawn using the igraph package in R \cite
{igraph} with the Fruchterman and Reingold layout \cite
{Fruchterman1991}). These are consistent with what we observed in
simulation studies: the Newman--Girvan modularity and the
degree-corrected block model likelihood perform better in a network
with high degree variation, and the Erdos--Renyi modularity is more
robust to degree variation than the block model likelihood.

All criteria were maximized by tabu search, but for modularities we
also computed the solutions based on the eigendecomposition of the
modularity matrix. Both solutions were worse that those found by tabu
search, but while for Newman--Girvan modularity the difference was
slight (the adjusted Rand index of 0.781 instead of 0.819),
eigendecomposition of the Erdos--Renyi modularity yielded a poor result
similar to that of block model likelihood (with adjusted Rand index
value of 0.092 instead of 0.819 by tabu search). This suggests that
Erdos--Renyi modularity is numerically less stable under high degree
variation, in addition to being theoretically not consistent. More
analysis of the eigendecomposition-based solutions is needed for both
types of modularities to understand conditions under which these
approximations work well.

%

\section{Summary and discussion}
\label{secsummary}

In this paper we developed a general tool for checking consistency of
community detection criteria under the degree-corrected stochastic
block model, a more general and practical model than the standard
stochastic block model for which such theory was previously available
\cite{Bickel&Chen2009}. This general tool allowed us to obtain
consistency results for four different community detection criteria,
and, to the best of our knowledge for the first time in the networks
literature, to clearly separate the effects of the model assumed for
criteria derivation from the model assumed true for analysis of the
criteria. What we have shown is, essentially, statistical common sense:
methods are consistent when the model they assume holds for the data.
The parameter constraints are needed when methods implicitly rely on
them, although we found that the two different modularity methods,
while using the same constraint in spirit, require somewhat different
conditions on parameters to be consistent. The theoretical analysis
agrees well with both simulation studies and the data analysis, which
also indicate that the methods with better theoretical consistency
properties do not always perform best in practice: there is a cost
associated with fitting the extra complexity of the degree-corrected
model, and if there is not enough data for that, or the data does not
have much variation in node degrees, simpler methods based on the
standard stochastic block model will in fact do better.

There are many questions that require further investigation here, even
in the context of model-based community detection when a model is
assumed true. For example, we assumed that $K$ is known, which is not
unreasonable in some cases (e.g., dividing political blogs into liberal
and conservative), but is in general a difficult open problem in
community detection. Standard methods such as AIC and BIC do not seem
to lend themselves easily to this case, because of parameters
disappearing in nonstandard ways when going from $K+1$ to $K$ blocks. A
permutation test was proposed in \cite{Zhaoetal2011}, but clearly more
work is needed.
There is also the question of what happens if $K$ is allowed to grow
with $n$, which is probably more realistic than fixed $K$; for the
stochastic block model, this case has been considered by \cite
{Choietal2011} and \cite{Rohe2011}, but their analysis is specific to
the particular methods they considered and does not extend easily to
the degree-corrected block model.
Another open question is the properties of approximate but more easily
computable solutions based on the eigendecomposition, as opposed to the
properties of global maximizers we studied here. For the stochastic
block model, part of this analysis was performed in \cite{Rohe2011}.
Our practical experience suggests that the behavior of eigenvectors can
be quite complicated, and it is not understood at this point when this
approximation works well. Finally, the sparse case $\lambda_n = O(1)$
is an open problem in general, although results for some special cases
of the stochastic block model have been recently obtained \cite
{coja09,Decelleetal2011}.

\begin{appendix}\label{app}

\section*{Appendix}
We start from summarizing notation. Let $R(\V{e}),V(\V{e}) \in
\mathcal
{R}^{K\times K \times M}$,
$\hat{\Pi} \in\mathcal{R}^{K \times M}$, $f(\V{e}),f^0(\V{e}) \in
\mathcal{R}^K$, where
\begin{eqnarray*}
R_{kau}(\V{e}) & =&\frac{1}{n}\sum_{i=1}^n
I(e_i=k,c_i=a,\theta_i=x_u) ,
\\
V_{kau}(\V{e}) & =&\frac{\sum_{i=1}^n I(e_i=k,c_i=a,\theta_i=x_u)}{\sum_{i=1}^n I(c_i=a,\theta_i=x_u)} ,
\\
\hat{\Pi}_{au} & =&\frac{1}{n}\sum_{i=1}^n
I(c_i=a,\theta_i=x_u) ,
\\
f_k(\V{e}) & =&\frac{1}{n}\sum_{i=1}^n
I(e_i=k)=\sum_{au} V_{kau}(\V
{e})\hat{\Pi}_{au} ,
\\
f^0_k(\V{e}) & =&\sum_{au}
V_{kau}(\V{e})\Pi_{au} .
\end{eqnarray*}
Even though the arbitrary labeling $\V{e}$ is not random, intuitively
one can think of~$R$ as the empirical joint distribution of $\V{e}$,
$\V
{c}$, and $\VV{\theta}$, $V$ as the conditional distribution of $\V{e}$
given $\V{c}$ and $\VV{\theta}$. Further, $\hat\Pi$ is the empirical
joint distribution of $\V{c}$ and $\VV{\theta}$, and thus an estimate of
their true joint distribution $\Pi$, $f$ is the empirical marginal
``distribution'' of $\V{e}$, and $f^0$ is the same marginal but with
the empirical joint distribution $\hat\Pi$ replaced by its population
version $\Pi$.
Then $\sum_{k} V_{kau}(\V{e})=1$, and $V_{kau}(\V{c}) = I(k=a)$ for all
$u$. Further, define $\hat{T}(\V{e}) \in\mathcal{R}^{K\times K}$ to be
a rescaled expectation of the matrix $O$ conditional on $c$ and $\theta$,
\[
\hat{T}_{kl}(\V{e})  =\frac{1}{\mu_n} \mathbb{E}
[O_{kl}| \V {c},\VV{\theta} ] .
\]
From Proposition \ref{prop1},
\begin{eqnarray*}
\hat{T}_{kl}(\V{e}) & = &\sum_{abuv}
x_u x_v P_{ab} R_{kau}(\V{e})
R_{lbv}(\V{e})
\\
& =& \sum_{abuv} x_u x_v
P_{st}V_{kau}(\V{e})\hat{\Pi}_{au} V_{lbv}(
\V {e})\hat{\Pi}_{bv} .
\end{eqnarray*}
Replacing $\hat{\Pi}$ by its expectation $\hat{\Pi}$, we define
$T(\V
{e}) \in\mathcal{R}^{K\times K}$ by
\[
T_{kl}(\V{e}) =\sum_{abuv}
x_u x_v P_{st}V_{kau}(\V{e})
\Pi_{au} V_{lbv}(\V{e})\Pi_{bv} .
\]
Also define $X(\V{e})\in\mathcal{R}^{K\times K} $ to be the rescaled
difference between $O$ and its conditional expectation,
\[
X_{kl}(\V{e}) =\frac{O_{kl}(\V{e})}{\mu_n}-\hat{T}_{kl}(\V{e}) .
\]
%
These quantities will be used in the proof of the general Theorem \ref
{genThm}, where we first approximate $\frac{1}{\mu_n} O_{kl}$ by
$\hat
{T}_{kl}(\V{e}) $ and then approximate $\hat{T}_{kl}(\V{e}) $ by
$T_{kl}(\V{e})$.

\begin{pf*}{Proof of Proposition \ref{prop1}}
We only proof \eqref{prop1a} since \eqref{prop1b} is trivial.
\begin{eqnarray*}
\hspace*{-4pt}&& \frac{1}{\mu_n} \mathbb{E} [O_{kl}| \V{c},\VV{\theta} ]
\\
\hspace*{-4pt}&&\qquad = \frac{1}{\mu_n} \sum
_{ij} \sum_{abuv} \mathbb {E}
\bigl[A_{ij}I(e_i=k,c_i=a,
\theta_i=x_u)I(e_j=l,c_j=b,
\theta_j=x_v)|\V {c},\VV{\theta}\bigr]
\\
\hspace*{-4pt}& &\qquad= \sum_{abuv} x_u x_v
P_{ab} R_{kau}(\V{e}) R_{lbv}(\V{e}) =
H_{kl}\bigl(R(\V{e})\bigr) .
\end{eqnarray*}
\upqed\end{pf*}

Before we proceed to the general theorem, we state a lemma based on
Bernstein's inequality. 
%
\begin{lemma}
\label{lemma-bernstein}
Let $\|X\|_{\infty}=\max_{kl} |X_{kl}|$ and $|\V{e}-\V{c}|=\sum_{i=1}^n
I(e_i\neq c_i)$. Then
%
\begin{equation}
\mathbb{P}\Bigl(\max_{\V{e}}\bigl\|X(\V{e})\bigr\|_{\infty} \geq\varepsilon
\Bigr) \leq 2K^{n+2} \exp \biggl( -\frac{1}{8C}
\varepsilon^2 \mu_n \biggr) \label {Bernstein1}
\end{equation}
for $ \varepsilon<3C$, where $C=\max\{x_u x_v P_{ab}\}$.
%
\begin{equation}\quad
\mathbb{P}\Bigl(\max_{|\V{e}-\V{c}|\leq m} \bigl\|X(\V{e})-X(\V{c})\bigr\|_{\infty
}\geq
\varepsilon\Bigr) \leq2 \pmatrix{{n}
\cr
{m}} K^{m+2} \exp \biggl( -
\frac
{3}{8} \varepsilon\mu_n \biggr) \label{Bernstein2}
\end{equation}
for $\varepsilon\geq6C m/n$.
%
\begin{equation}
\mathbb{P}\Bigl(\max_{|\V{e}-\V{c}|\leq m} \bigl\|X(\V{e})-X(\V{c})\bigr\|_{\infty
}\geq
\varepsilon\Bigr) \leq2 \pmatrix{{n}
\cr
{m}} K^{m+2} \exp \biggl( -
\frac
{n}{16 m C}\varepsilon^2 \mu_n \biggr)\hspace*{-35pt}
\label{Bernstein3}
\end{equation}
for $\varepsilon< 6C m/n$.
\end{lemma}
This lemma is similar to Lemma 1.1 of \cite{Bickel&Chen2009}, with a
few minor errors corrected. The proof can be found in the electronic
supplement to this article \cite{Zhaoetal-aossupp}.
\begin{pf*}{Proof of Theorem \ref{genThm}}
The proof is divided into three steps.

\textit{Step} 1: show that $F  (\frac{O(\V{e})}{\mu_n},f(\V{e})
)$ is uniformly close to its population version. More precisely,
we need to prove that there exists $\varepsilon_n \rightarrow0$, such that
%
\begin{equation}
\label{step1}\qquad \mathbb{P} \biggl(\max_{\V{e}} \biggl\llvert F \biggl(
\frac{O(\V{e})}{\mu_n},f(\V{e}) \biggr)-F\bigl(T(\V{e}),f^0(\V{e})\bigr)
\biggr\rrvert < \varepsilon_n \biggr) \rightarrow1 \qquad \mbox{if }
\lambda_n \rightarrow\infty.
\end{equation}

Since
\begin{eqnarray*}
\biggl| F \biggl(\frac{O(\V{e})}{\mu_n},f(\V{e}) \biggr) -  F\bigl(T(\V
{e}),f^0(\V {e})\bigr) \biggr| &\leq&\biggl| F \biggl(\frac{O(\V{e})}{\mu_n},f(\V{e})
\biggr) -F\bigl(\hat{T}(\V{e}),f(\V{e})\bigr) \biggr|
\\
& &{} + \bigl|F\bigl(\hat{T}(\V{e}),f(\V{e})\bigr)-F\bigl(T(\V{e}),f^0(\V{e})
\bigr) \bigr| ,
\end{eqnarray*}
it is sufficient to bound these two terms uniformly. By Lipschitz continuity,
%
\begin{equation}
\biggl\llvert F \biggl(\frac{O(\V{e})}{\mu_n},f(\V{e}) \biggr)-F\bigl(\hat {T}(\V
{e}),f(\V{e})\bigr) \biggr\rrvert \leq M_1 \bigl\|X(e)\bigr\|_{\infty} .
\label{lip}
\end{equation}
By \eqref{Bernstein1}, \eqref{lip} converges to 0 uniformly if
$\lambda_n \rightarrow\infty$, and
%
\begin{eqnarray}\label{lip2}
&&\bigl|F\bigl(\hat{T}(\V{e}),f(\V{e})\bigr)-F\bigl(T(\V{e}),f^0(\V{e})
\bigr) \bigr|
\nonumber
\\[-8pt]
\\[-8pt]
\nonumber
&&\qquad \leq M_1 \bigl\|\hat{T}(\V{e})-T(\V{e})\bigr\|_{\infty} +
M_2 \bigl\|f(\V {e})-f^0(\V {e})\bigr\|
\end{eqnarray}
where $\| \cdot\|$ is the Euclidean norm for vectors.
Further,
%
\begin{eqnarray}
\bigl|\hat{T}_{kl}(\V{e})-T_{kl}(\V{e})\bigr|& = &\biggl |\sum
_{abuv} x_u x_v P_{ab}
V_{kau}(\V{e}) V_{lbv}(\V{e}) (\hat{\Pi}_{au}\hat{
\Pi}_{bv}-\Pi_{au}\Pi_{bv}) \biggr|
\nonumber
\\[-8pt]
\\[-8pt]
\nonumber
&\leq& \sum_{abuv} x_u x_v
P_{ab} |\hat{\Pi}_{au}\hat{\Pi }_{bv}-
\Pi_{au}\Pi_{bv} | ,
\end{eqnarray}
and
%
\begin{equation}
\bigl|f_k(\V{e})-f^0_k(\V{e})\bigr| = \biggl\llvert
\sum_{au} V_{kau}(\V{e}) (\hat {\Pi
}_{au}-\Pi_{au}) \biggr\rrvert \leq\sum
_{au} |\hat{\Pi}_{au}-\Pi_{au}| .
\end{equation}
Since $\hat{\Pi} \stackrel{P}{\rightarrow} \Pi$, \eqref{lip2} converges
to 0 uniformly. Thus, \eqref{step1} holds.

\textit{Step} 2: Prove that there exists $\delta_n\rightarrow0$,
such that
%
\begin{equation}
\label{step2} \mathbb{P} \biggl( \max_{\{\V{e}\dvtx \|V(\V{e})-\mathbb{I}\|_1\geq
\delta_n\}
} F \biggl(
\frac{O(\V{e})}{\mu_n},f(\V{e}) \biggr) < F \biggl( \frac{O(\V
{c})}{\mu_n},f(\V{c})
\biggr) \biggr) \rightarrow1 ,
\end{equation}
where $\|W\|_1=\sum_{kau} |W_{kau}|$ for $W \in\mathcal{R}^{K\times K
\times M}$.

By continuity and ($\ast$), there exists $\delta_n\rightarrow0$,
such that
\[
F\bigl(T(\V{c}),f^0(\V{c})\bigr)-F\bigl(T(\V{e}),f^0(
\V{e})\bigr)>2\varepsilon_n
\]
if $\|V(\V{e})-\mathbb{I}\|_1\geq\delta_n$, where $\mathbb{I}=V(\V
{c})$. Thus, from \eqref{step1},

\begin{eqnarray*}
&&\mathbb{P} \biggl( \max_{ \{\V{e}\dvtx  \|V(\V{e})-\mathbb{I}\|_1 \geq
\delta
_n \} }  F \biggl( \frac{O(\V{e})}{\mu_n},f(\V{e})
\biggr) < F \biggl( \frac
{O(\V
{c})}{\mu_n},f(\V{c}) \biggr) \biggr)
\\
&&\qquad\geq\mathbb{P} \biggl( \biggl| \max_{ \{\V{e}\dvtx  \|V(\V{e})-\mathbb{I}\|_1
\geq\delta_n \} }  F\biggl( \frac
{O(\V
{e})}{\mu_n},f(\V{e})
\biggr)\\
&&\hspace*{27pt}\qquad{}-\max_{ \{\V{e}\dvtx  \|V(\V{e})-\mathbb
{I}\|_1
\geq\delta_n \} } F\bigl(T(\V{e}),f^0(\V{e})\bigr) \biggr| <
\varepsilon_n,
\\
&&\hspace*{46pt}\biggl | F \biggl( \frac{O(\V{c})}{\mu_n},f(\V{c}) \biggr)-F\bigl(T(\V
{c}),f^0(\V{c})\bigr) \biggr| < \varepsilon_n \biggr) \rightarrow1.
\end{eqnarray*}
\eqref{step2} implies
\begin{eqnarray*}
\mathbb{P} \bigl( \bigl\Vert V(\hat{c})-\mathbb{I}\bigr\Vert < \delta_n \bigr)
\rightarrow1 .
\end{eqnarray*}
Since
\begin{eqnarray*}
\frac{1}{n}|\V{e}-\V{c}| & =& \frac{1}{n} \sum
_{i=1}^n I(c_i \neq e_i) =
\sum_{au} \Pi_{au}\bigl(1-V_{aau}(
\V{e})\bigr) 
\leq\sum_{au}
\bigl(1-V_{aau}(\V{e})\bigr)
\\
& =& \frac{1}{2} \biggl(\sum_{au}
\bigl(1-V_{aau}(\V{e})\bigr)+\sum_{au} \sum
_{k\neq a} V_{kau}(\V{e}) \biggr) =
\frac{1}{2} \bigl\|V(\V{e})-\mathbb{I}\bigr\|_1 ,
\end{eqnarray*}
weak consistency follows.

\textit{Step} 3: In order to prove strong consistency, we need to show that
%
\begin{equation}
\mathbb{P} \biggl( \max_{ \{\V{e}\dvtx 0<\|V(\V{e})-\mathbb{I}\|_1 <
\delta_n
\}} F \biggl( \frac{O(\V{e})}{\mu_n},f(\V{e})
\biggr) < F \biggl( \frac
{O(\V{c})}{\mu_n},f(\V{c}) \biggr) \biggr) \rightarrow1.\hspace*{-35pt}
\label{step3}
\end{equation}

Note that combining \eqref{step2} and \eqref{step3}, we have
\[
\mathbb{P} \biggl( \max_{ \{\V{e}\dvtx  \V{e} \neq\V{c} \}} F \biggl( \frac
{O(\V{e})}{\mu_n},f(\V{e})
\biggr) < F \biggl( \frac{O(\V
{c})}{\mu_n},f(\V{c}) \biggr) \biggr) \rightarrow1,
\]
which implies the strong consistency.

Here we closely follow the derivation given in \cite{AiroldiNotes}.
To prove \eqref{step3}, note that by Lipschitz continuity and the
continuity of derivatives of $F$ with respect to $V(\V{e})$ in the
neighborhood of $\mathbb{I}$, we have
%
\begin{eqnarray}\label{lip3}
&& F \biggl( \frac{O(\V{e})}{\mu_n},f(\V{e}) \biggr)-F \biggl( \frac
{O(\V
{c})}{\mu_n},f(
\V{c}) \biggr)
\nonumber
\\[-8pt]
\\[-8pt]
\nonumber
&&\qquad=F\bigl(\hat{T}(\V{e}),f(\V{e})\bigr)-F\bigl(\hat {T}(\V {c}),f(
\V{c})\bigr)+\Delta(\V{e},\V{c}),
\end{eqnarray}
where $|\Delta(\V{e},\V{c})| \leq M'(\|X(\V{e})- X(\V{c})\|_{\infty})$, and
%
\begin{eqnarray}
&&F\bigl(T(\V{e}),f^0(\V{e})\bigr)-F\bigl(T(\V{c}),f^0(
\V{c})\bigr)
\nonumber
\\[-8pt]
\\[-8pt]
\nonumber
&&\qquad\leq-C'\bigl\|V(\V {e})-\mathbb {I}\bigr\|_1+o\bigl(
\bigl\|V(\V{e})-\mathbb{I}\bigr\|_1\bigr).
\end{eqnarray}
Since the derivative of $F$ is continuous with respect to $V(\V{e})$ in
the neighborhood of $\mathbb{I}$, there exists a $\delta'$ such that
%
\begin{eqnarray}\label{taylor}
&& F\bigl(\hat{T}(\V{e}),f(\V{e})\bigr)-F\bigl(\hat{T}(\V{c}),f(\V{c})\bigr)
\nonumber
\\[-8pt]
\\[-8pt]
\nonumber
&&\qquad \leq -
\bigl(C'/2\bigr)\bigl\|V(\V {e})-\mathbb{I}\bigr\|_1+o\bigl(\bigl\|V(
\V{e})-\mathbb{I}\bigr\|_1\bigr)
\end{eqnarray}
holds when $\|\hat{\Pi}-\Pi\|_{\infty}\leq\delta'$. Since $\hat
{\Pi}
\rightarrow\Pi$, \eqref{taylor} holds with probability approaching 1.
Combining \eqref{lip3} and \eqref{taylor}, it is easy to see that
\eqref
{step3} follows if we can show
%
\begin{equation}
\mathbb{P}\Bigl(\max_{\{\V{e}\neq\V{c}\}}\bigl|\Delta(\V{e},\V{c})\bigr|\leq C'\bigl\|V(
\V {e})-\mathbb{I}\bigr\|_1 /4\Bigr) \rightarrow1.
\end{equation}
Again note that $\frac{1}{n}|\V{e}-\V{c}| \leq\frac{1}{2} \|V(\V
{e})-\mathbb{I}\|_1$.
So for each $m \geq1$,
%
\begin{eqnarray}\label{hhh}
&&\mathbb{P}\Bigl(\max_{|\V{e}-\V{c}|=m} \bigl|\Delta(\V{e},\V{c})\bigr|>C'\bigl\|V(\V
{e})-\mathbb{I}\bigl\|_1 /4\Bigr)
\nonumber
\\[-8pt]
\\[-8pt]
\nonumber
&&\qquad \leq\mathbb{P}\biggl(
\max_{|\V{e}-\V{c}|\leq m} \bigl\|X(\V{e})-X(\V{c})\bigr\|_{\infty}>\frac{C'm}{2M'n}
\biggr)=I_1 .
\end{eqnarray}
Let $\alpha=C'/2M'$, if $\alpha\geq6C$, by \eqref{Bernstein2},
\begin{eqnarray*}
I_1 & \leq&2 K^{m+2} n^m \exp \biggl( -\alpha
\frac{3m}{8n} \mu_n \biggr)
\\
& =& 2 K^2 \bigl[K \exp\bigl(\log n - \alpha\mu_n/(8/3 n)
\bigr) \bigr]^m.
\end{eqnarray*}
If $\alpha<6C$, by \eqref{Bernstein3},
\begin{eqnarray*}
I_1 & \leq& 2 K^{m+2} n^m \exp \biggl( -
\alpha^2 \frac{m}{16C n} \mu_n \biggr)
\\
& =& 2 K^2 \bigl[K \exp\bigl(\log n - \alpha^2
\mu_n /(16C n)\bigr) \bigr]^m.
\end{eqnarray*}

In both cases, since $\lambda_n/\log n \rightarrow\infty$,
\begin{eqnarray*}
&&\mathbb{P}\Bigl(\max_{\{\V{e}\neq\V{c}\}}\bigl|\Delta(\V{e},\V{c})\bigr|> C'\bigl\| V(\V
{e})-\mathbb{I}\bigr\|_1 /4\Bigr)
\\
&&\qquad= \sum_{m=1}^{\infty}
\mathbb{P}\Bigl(\max_{|\V
{e}-\V
{c}|=m} \bigl|\Delta(\V{e},\V{c})\bigr|>C'\bigl\|V(
\V{e})-\mathbb{I}\bigr\|_1 /4\Bigr) \rightarrow0
\end{eqnarray*}
as $n\rightarrow\infty$, which completes the proof.
\end{pf*}

\begin{pf*}{Proof of Theorem \ref{corERM}}
The regularity conditions are easy to verify. To check the key
condition ($\ast$), note that under the block model assumption, ($\ast
$) becomes
\begin{itemize}[($\ast\ast$)]
\item[($\ast\ast$)] $F(H(S),h(S))$ is uniquely maximized over
$\mathscr
{S}=\{S\dvtx  S \geq0, \sum_k S_{ka}=\pi_{a} \}$ by $S=D$, with
$D=\operatorname{diag}(\VV{\pi})$,
\end{itemize}
where $S$ is a generic $K$ by $K$ matrix.

Up to a constant, the population version of $Q_{\mathrm{ERM}} $ is
\[
F\bigl(H(S),h(S)\bigr)=\sum_{k}
\bigl(H_{kk}-h_k^2 P_0\bigr).
\]
Using the identity,
\[
\sum_{k} \bigl(H_{kk}-h_k^2
P_0\bigr) +\sum_{k\neq l}
(H_{kl}-h_kh_l P_0)=\sum
_{kl} H_{kl} -\biggl(\sum
_{k} h_k\biggr)^2 P_0
=0,
\]
and define
\begin{eqnarray*}
\Delta_{kl}= \cases{ %
1, &\quad $\mbox{if $k=l$}$,
\vspace*{2pt}
\cr
-1, &\quad $\mbox{if $k \neq l$.}$}
\end{eqnarray*}
Then we have
\begin{eqnarray*}
F\bigl(H(S),h(S)\bigr) & =&\frac{1}{2} \sum_{kl}
\Delta_{kl} (H_{kl}-h_k h_l
P_0)\\
& =&\frac{1}{2} \sum_{kl}
\Delta_{kl} \biggl(\sum_{ab}
S_{ka} S_{lb} P_{ab} -\sum
_{ab} S_{ka} S_{lb} P_0\biggr)
\\
& =&\frac{1}{2} \sum_{kl} \sum
_{ab} S_{ka} S_{lb}\Delta_{kl} (
P_{ab}-P_0)\\
& \leq&\frac{1}{2} \sum
_{kl} \sum_{ab} S_{ka}
S_{lb} \Delta_{ab} (P_{ab}-P_0)
\\
& = &\frac{1}{2} \sum_{ab}
\Delta_{ab} \pi_a \pi_b (P_{ab}-P_0)
= F\bigl(H(D),h(D)\bigr).
\end{eqnarray*}
Now it remains to show the diagonal matrix $D$ (up to a permutation) is
the unique maximizer of $F$. This follows from Lemma 3.2 in \cite
{Bickel&Chen2009}, since equality holds only if $\Delta_{kl}=\Delta_{ab}$ when
$S_{ka}S_{lb}>0$ and $\Delta$ does not have two identical columns.
\end{pf*}

\begin{pf*}{Proof of Theorem \ref{corNGM}}
The consistency of Newman--Girvan modularity under the block model has
already been shown in \cite{Bickel&Chen2009}. To extend this result to
the degree-corrected block model, define $\tilde{S}_{ka}=\sum_u x_u
S_{kau}$. Then
\begin{eqnarray*}
\tilde{\pi}_a & =&\sum_k
\tilde{S}_{ka},
\\
H_{kl} & =& \sum_{abuv} x_u
x_v P_{ab}S_{kau}S_{lbv} = \sum
_{ab} \tilde {S}_{ka} \tilde{S}_{lb}
P_{ab},
\\
H_{k} & =& \sum_{l} H_{kl} =
\sum_{as} \tilde{S}_{ka} \tilde{\pi
}_{s} P_{as}.
\end{eqnarray*}
The population version of $Q_{\mathrm{NGM}}$ is
\[
F\bigl(H(S)\bigr)  = \sum_{k} \biggl(
\frac{H_{kk}}{\tilde{P}_0}- \biggl(\frac
{H_k}{\tilde{P}_0} \biggr)^2 \biggr).
\]
Using the identity
\[
\sum_{k} \biggl( \frac{H_{kk}}{\tilde{P}_0}- \biggl(
\frac
{H_k}{\tilde
{P}_0} \biggr)^2 \biggr) +\sum
_{k\neq l} \biggl( \frac
{H_{kl}}{\tilde
{P}_0}- \frac{H_kH_l}{\tilde{P}_0^2} \biggr)=
\sum_{kl} \frac
{H_{kl}}{\tilde{P}_0} - \biggl(\sum
_{k} \frac{H_k}{\tilde{P}_0} \biggr)^2 =0,
\]
we obtain
\begin{eqnarray*}
F\bigl(H(S)\bigr) & =& \frac{1}{2} \sum_{kl}
\Delta_{kl} \biggl( \frac{\sum_{ab}
\tilde{S}_{ka}\tilde{S}_{lb}P_{ab}}{\tilde{P}_0}-\frac{(\sum_{as}\tilde
{S}_{ka}\tilde{\pi}_s P_{as})(\sum_{bt}\tilde{S}_{lb}\tilde{\pi}_t
P_{bt})}{\tilde{P}_0^2} \biggr)
\\
& =&\frac{1}{2} \sum_{kl} \sum
_{ab} \tilde{S}_{ka}\tilde{S}_{lb}
\Delta_{kl} \biggl( \frac{P_{ab}}{\tilde{P}_0}- \frac{(\sum_{s}\tilde
{\pi}_s
P_{as})(\sum_{t}\tilde{\pi}_t P_{bt})}{\tilde{P}_0^2} \biggr)
\\
& \leq&\frac{1}{2} \sum_{kl} \sum
_{ab} \tilde{S}_{ka}\tilde{S}_{lb}
\Delta_{ab} \biggl( \frac{P_{ab}}{\tilde{P}_0}- \frac{(\sum_{s}\tilde
{\pi}_s P_{as})(\sum_{t}\tilde{\pi}_t P_{bt})}{\tilde
{P}_0^2} \biggr)
\\
& =& \frac{1}{2} \sum_{ab}
\Delta_{ab} \tilde{\pi}_a\tilde{\pi}_b \biggl(
\frac{P_{ab}}{\tilde{P}_0}- \frac{(\sum_{s}\tilde{\pi}_s
P_{as})(\sum_{t}\tilde{\pi}_t P_{bt})}{\tilde{P}_0^2} \biggr) =F\bigl(H(\mathbb{D})\bigr).
\end{eqnarray*}
Similar to Theorem \ref{corERM}, $D$ is the unique maximizer of
$F(H(\tilde{S}))$, so it is enough to show $S=\mathbb{D}$ whenever
$\tilde{S}=D$ to prove uniqueness. $\tilde{S}=D$ implies $\tilde
{S}_{ka}=0$, if $k\neq a$. Since $x_u>0$, we obtain $S_{kau}=0$ if
$k\neq a$, which gives the result.

We note that this argument cannot be applied to prove the consistency
of Erdos--Renyi modularity under degree-corrected block models, because
in that case $h_k=\sum_{au} S_{kau} \neq\sum_{a}(\sum_u
x_u S_{kau}) = \sum_{a} \tilde{S}_{ka}$, when we use the transformation
$\tilde{S}_{ka}=\sum_u x_u S_{kau}$.
\end{pf*}

\begin{pf*}{Proof of Theorem \ref{corBL}}
Up to a constant, the population version of $Q_{\mathrm{BL}}$ is
\[
F\bigl(H(S),h(S)\bigr)=\sum_{kl} \biggl(
H_{kl} \log\frac{H_{kl}}{h_kh_l}-H_{kl} \biggr).
\]
Let $g_{kl}=H_{kl}/(h_kh_l)$,
\begin{eqnarray*}
F\bigl(H(S),h(S)\bigr) &=&\sum_{kl}
(H_{kl} \log g_{kl}-h_kh_l
g_{kl}) 
=\sum_{abkl}
S_{ka}S_{lb} (P_{ab} \log g_{kl}-
g_{kl})
\\
& \leq&\sum_{ab}\sum_{kl}
S_{ka}S_{lb}(P_{ab} \log P_{ab} -
P_{ab})
\\
& =& \sum_{ab} (\pi_a
\pi_b P_{ab} \log P_{ab}- \pi_a
\pi_b P_{ab})=F\bigl(H(D),h(D)\bigr).
\end{eqnarray*}
%
Since the inequality holds if and only if $g_{kl}=P_{ab}$ when
$S_{ka}S_{lb}>0$, uniqueness follows from Lemma \ref{lemma-uniq},
stated next.
\end{pf*}
%
\begin{lemma}
\label{lemma-uniq}
Let $g$, $P$, $S$ be $K \times K$ matrices with nonnegative entries.
Assume that:
\begin{longlist}[(a)]
\item[(a)] $P$ and $g$ are symmetric;
\item[(b)] $P$ does not have two identical columns;
\item[(c)] there exists at least one nonzero entry in each column of $S$;
\item[(d)] for $1 \leq k,l,a,b \leq K, g_{kl}=P_{ab}$ whenever
$S_{ka}S_{lb}>0$.
\end{longlist}
Then $S$ is a diagonal matrix or a row/column permutation of a diagonal matrix.
\end{lemma}
This lemma is a generalization of Lemma 3.2 in \cite{Bickel&Chen2009}.
The proof is given in the electronic supplement \cite{Zhaoetal-aossupp}.
%
%
%
\begin{pf*}{Proof of Theorem \ref{corDCBL}}
Up to a constant, the population version of $Q_{\mathrm{DCBM}}$ is
%
\begin{equation}
F\bigl(H(S)\bigr)=\sum_{kl} \biggl(
H_{kl} \log\frac{H_{kl}}{H_kH_l}-H_{kl} \biggr),
\end{equation}
where we only check ($\ast\ast$) [the form ($\ast$) takes under the block
model]. The generalization to the degree-corrected block model is
similar to the proof of Theorem~\ref{corNGM} and is omitted.

Let $g_{kl}=H_{kl}/(H_kH_l)$, and
\begin{eqnarray*}
F\bigl(H(S)\bigr) & =& \sum_{kl} (H_{kl}
\log g_{kl}-H_kH_l g_{kl})
\\
& =& \sum_{kl} \biggl[\sum
_{ab} S_{ka}S_{lb}P_{ab} \log
g_{kl}-\biggl(\sum_{as}S_{ka}
\pi_{s} P_{as} \biggr) \biggl(\sum
_{bt} \pi_t S_{lb} P_{tb}
\biggr) g_{kl}
\biggr]
\\
& =& \sum_{kl} \sum_{ab}
S_{ka}S_{lb} \biggl[P_{ab} \log g_{kl}-
\biggl(\sum_{s} \pi_{s} P_{as}
\biggr) \biggl(\sum_{t} \pi_t
P_{tb}\biggr) g_{kl} \biggr] =I_2 .
\end{eqnarray*}
Since $\operatorname{arg\,max}_{x} (c_1\log x-c_2 x)=c_1/c_2$, replacing
$g_{kl}$ by
\[
\frac{P_{ab}}{(\sum_{s} \pi_{s} P_{as} )(\sum_{t} \pi_t P_{tb})},
\]
we obtain
\begin{eqnarray*}
I_2 &\leq& \sum_{kl}\sum
_{ab} S_{ka}S_{lb} \biggl[
P_{ab} \log\frac
{P_{ab}}{(\sum_{s} \pi_{s} P_{as} )(\sum_{t} \pi_t P_{tb})}- P_{ab} \biggr]
\\
&= & \sum_{ab} \biggl[\pi_a
\pi_b P_{ab} \log\frac{P_{ab}}{(\sum_{s} \pi_{s} P_{as} )(\sum_{t} \pi_t P_{tb})}-\pi_a
\pi_b P_{ab} \biggr] = F\bigl(H(D)\bigr).
\end{eqnarray*}
\upqed\end{pf*}
\end{appendix}

\section*{Acknowledgments}
This work was carried out while Yunpeng Zhao was a Ph.D. student at the
University of Michigan. We thank Brian Karrer and Mark Newman
(University of Michigan) for helpful comments and corrections, Peter
B\"uhlmann (ETH) for the role he played as Editor, and two anonymous
referees for their constructive feedback and corrections.

\begin{supplement}
\stitle{Proofs of Lemmas \ref{lemma-bernstein} and \ref{lemma-uniq}.}
\slink[doi]{10.1214/12-AOS1036SUPP} 
\sdatatype{.pdf}
\sfilename{aos1036\_supp.pdf}
\sdescription{The supplemental material contains proofs of Lemmas \ref
{lemma-bernstein} and \ref{lemma-uniq} stated in the \hyperref
[app]{Appendix}.}
\end{supplement}

%
%

%

\printaddresses

\end{document}


\begin{frontmatter}

\title{Correction to the proof of consistency of community detection}
\runtitle{Correction to consistency of community detection}

\begin{aug}
\author[A]{\fnms{Peter~J.}~\snm{Bickel}\ead[label=pb]{bickel@stat.berkeley.edu}},
\author[B]{\fnms{Aiyou}~\snm{Chen}\ead[label=ac]{aiyouchen@google.com}},
\author[C]{\fnms{Yunpeng}~\snm{Zhao}\ead[label=yz]{yzhao15@gmu.edu}},
\author[D]{\fnms{Elizaveta}~\snm{Levina}\corref{}\ead[label=el]{elevina@umich.edu}}
\and
\author[D]{\fnms{Ji}~\snm{Zhu}\ead[label=jz]{jizhu@umich.edu}}
\runauthor{P.~J. Bickel et al.}
\affiliation{University of California, Berkeley,
Google Inc, George Mason University, University of Michigan and
University of Michigan}

\address[A]{P.~J. Bickel\\
Department of Statistics\\
University of California, Berkeley\\
367 Evans Hall\\
Berkeley, California 94720-3860\\
USA\\
\printead{pb}}

\address[B]{A. Chen\\
Google Inc\\
1600 Amphitheatre Pkwy\\
Mountain View, California 94043\\
USA\\
\printead{ac}}

\address[C]{Y. Zhao\\
Department of Statistics\\
George Mason University\\
1714 Engineering Building \\
4400 University Drive\\
Fairfax, Virginia 22030-4444\\
USA\\
\printead{yz}}

\address[D]{E. Levina\\
J. Zhu\\
Department of Statistics\\
University of Michigan\\
311 West Hall\\
1085 S. University Ave.\\
Ann Arbor, Michigan 48109-1107\hspace*{3pt}\\
USA\\
\printead{el}\\
\phantom{E-mail:\ }\printead*{jz}}

\end{aug}

\received{\smonth{8} \syear{2014}}
\revised{\smonth{9} \syear{2014}}

%
\begin{abstract}
This note corrects an error in two related proofs of consistency of community
detection: under stochastic block models by Bickel and Chen
[\textit{Proc. Natl. Acad. Sci. USA} \textbf{106} (2009)
21068--21073]
and under degree-corrected stochastic block model by Zhao, Levina and
Zhu [\textit{Ann. Statist.} \textbf{40} (2012) 2266--2292].
\end{abstract}

\begin{keyword}[class=AMS]
\kwd{62G20}
\end{keyword}
\begin{keyword}
\kwd{Network communities}
\kwd{stochastic block model}
\kwd{degree-corrected stochastic block model}
\kwd{consistency of community detection}
\end{keyword}
\end{frontmatter}
\setcounter{section}{1}

This note provides a correction to the proof of consistency of
community detection under
degree-corrected stochastic block models \cite{Zhaoetal2012},
published in this journal. The same error appeared earlier in the
proof of consistency under the stochastic block models
\cite{BickelChen2009}. In this note, we provide the correction for
the proof of \cite{Zhaoetal2012}, using the notation of that paper,
since the case of the
degree-corrected stochastic block models is more general and includes
the regular stochastic block models as a special case. Very similar
arguments can be used to correct the proof of \cite{BickelChen2009}
directly.

We start by very briefly restating notation. Let $\mathbf{e}$ be an
arbitrary set of label assignments, $\mathbf{c}$ be the true label
assignments and $\hat{\mathbf{c}}$ be the maximizer of a community detection
criterion. Let $O(\mathbf{e}) \in\mathcal{R}^{K\times K}$,
$V(\mathbf{e}) \in
\mathcal{R}^{K\times K \times M}$, $\hat{\Pi} \in\mathcal{R}^{K
\times
M}$, $f(\mathbf{e}) \in\mathcal{R}^K$, where
\begin{eqnarray*}
O_{kl}(\mathbf{e}) & =&\sum_{ij}
A_{ij}I\{e_i=k,e_j=l\},
\\
V_{kau}(\mathbf{e}) & =&\frac{\sum_{i=1}^n I(e_i=k,c_i=a,\theta
_i=x_u)}{\sum_{i=1}^n I(c_i=a,\theta_i=x_u)},
\\
\hat{\Pi}_{au} & =&\frac{1}{n}\sum_{i=1}^n
I(c_i=a,\theta_i=x_u),
\\
f_k(\mathbf{e}) & =&\frac{1}{n}\sum
_{i=1}^n I(e_i=k)=\sum
_{au} V_{kau}(\mathbf{e})\hat{\Pi}_{au}.
\end{eqnarray*}
We considered community detection criteria that can be written in the
form
\[
Q(\mathbf{e})=F \biggl( \frac{O(\mathbf{e})}{\mu_n},f(\mathbf{e}) \biggr),
\]
where $\mu_n=n^2 \rho_n$ and $\rho_n \rightarrow0$ is the average
probability of an
edge in the network. For any matrix $B$,
$\|B\|_{\infty}=\max_{kl} |B_{kl}|$.

The statement $|\Delta(\mathbf{e},\mathbf{c})| \leq M_1(\|X(\mathbf{e})-
X(\mathbf{c})\|_{\infty})$ below (A.11) in \cite{Zhaoetal2012} is
incorrect. (We have replaced $M'$ and $C'$ in the original with
$M_1$ and $C_1$ in this correction since we will need more
constants.) For the proof to go through, we need a different way of
proving
%
\begin{equation}
\label{claim} \mathbb{P} \Bigl(\max_{1\leq|\mathbf{e}-\mathbf{c}| \leq\delta_n
n}\bigl|\Delta (\mathbf{e},
\mathbf{c})\bigr|-C_1\bigl\|V(\mathbf{e})-\mathbb{I}\bigr\|_1 /4\leq0
\Bigr) \rightarrow 1,
\end{equation}
where $\delta_n \rightarrow0$. Note that \eqref{claim} is similar to
the (A.14) in \cite{Zhaoetal2012}, with an extra constraint
$|\mathbf{e}-\mathbf{c}| \leq\delta_n n$. Since we have already proved
$\mathbb{P}(\frac{1}{n}|\hat{\mathbf{c}}-\mathbf{c}|\leq\delta
_n)\rightarrow1$ in
\cite{Zhaoetal2012}, \eqref{claim} will complete the proof, and the
conclusion of Theorem~4.1 in \cite{Zhaoetal2012} remains valid.

We first need a lemma based on Bernstein's inequality.

\begin{lemma}
\label{lemma-bernstein}
For $m \in\{1,\ldots,n\}$,
%
\begin{equation}
\label{Bernstein1} \mathbb{P}\Bigl(\max_{|\mathbf{e}-\mathbf{c}|\leq m}
\bigl\|X(\mathbf{e})\bigr\|
_{\infty} \geq \varepsilon \Bigr) \leq2 \pmatrix{n
\cr
m} K^{m+2}
\exp \biggl( -\frac{3 \mu_n
\varepsilon
^2}{4(\varepsilon+3)} \biggr).
\end{equation}
\end{lemma}

The proof of Lemma~\ref{lemma-bernstein} closely follows the proof of
(A.2) and (A.3) in \cite{Zhaoetal2012} and hence is omitted here.

Proof of \eqref{claim}:

By Taylor's expansion,
\begin{eqnarray*}
&&F \biggl( \frac{O(\mathbf{e})}{\mu_n},f(\mathbf{e}) \biggr)- F\bigl(\hat{T}(\mathbf{e}),f(
\mathbf{e})\bigr)\\
&&\qquad = \frac{\partial F}{\partial M}
\bigg|_{M=\hat{T}(\mathbf{e}), \mathbf{t}=f(\mathbf{e})} \operatorname{vec}\bigl(X(\mathbf{e})
\bigr)+O\bigl(\bigl\| X(\mathbf{e}) \bigr\|^2_\infty\bigr),
\end{eqnarray*}
where $\frac{\partial F}{\partial M}$ is the partial derivative over
the first component (vectorized) of $F(M, \mathbf{t})$.
Similarly,
\begin{eqnarray*}
&&F \biggl( \frac{O(\mathbf{c})}{\mu_n},f(\mathbf{c}) \biggr)- F\bigl(\hat {T}(
\mathbf{c}),f(\mathbf{c})\bigr)\\
&&\qquad = \frac{\partial F}{\partial M} \bigg|_{M=\hat{T}(\mathbf{c}), \mathbf{t}=f(\mathbf{c}) }\operatorname{vec}\bigl(X(
\mathbf {c})\bigr)+O\bigl(\bigl\|X(\mathbf{c}) \bigr\|^2_\infty\bigr).
\end{eqnarray*}
Since $\frac{\partial F}{\partial M}$ is continuous with respect to $M$
and $t$, and $ \hat{T}(\mathbf{e})$ and $f(\mathbf{e})$ are
continuous with
respect to $\mathbf{e}$, %
\begin{equation}
\frac{\partial F}{\partial M} \bigg|_{M=\hat{T}(\mathbf{e}),
\mathbf{t}=f(\mathbf{e})
}=\frac{\partial F}{\partial M} \bigg|_{M=\hat{T}(\mathbf{c}),
\mathbf{t}=f(\mathbf{c}) }+O\bigl(
\bigl\|V(\mathbf{e})-\mathbb{I}\bigr\|_1\bigr).
\end{equation}
Therefore, since
\begin{eqnarray*}
\Delta(\mathbf{e},\mathbf{c}) & =& F \biggl( \frac{O(\mathbf
{e})}{\mu_n},f(\mathbf{e})
\biggr)- F\bigl(\hat{T}(\mathbf{e}),f(\mathbf{e})\bigr) - F \biggl(
\frac{O(\mathbf{c})}{\mu_n},f(\mathbf{c}) \biggr) + F\bigl(\hat {T}(\mathbf{c}),f(
\mathbf{c})\bigr)
\nonumber
\\
&= & \frac{\partial F}{\partial M} \bigg|_{M=\hat{T}(\mathbf{c}),
\mathbf{t}=f(\mathbf{c}) }\operatorname{vec}\bigl(X(\mathbf{e})-X(\mathbf{c})\bigr)
+O\bigl(\bigl\| V(\mathbf{e})-\mathbb{I}\bigr\| _1\bigr)
\operatorname{vec}\bigl(X(\mathbf{e})
\bigr)
\nonumber
\\
&&{}+O\bigl(\bigl\|X(\mathbf{e}) \bigr\|^2_\infty\bigr) +O\bigl(\bigl\|X(
\mathbf{c})\bigr \|^2_\infty\bigr),
\end{eqnarray*}
we have
\begin{eqnarray*}
\bigl| \Delta(\mathbf{e},\mathbf{c})\bigr |& \leq &M_1\bigl\|X(\mathbf{e})- X(
\mathbf{c})\bigr\|_{\infty
}+M_2\bigl\| V(\mathbf{e})-\mathbb{I}
\bigr\|_1 \bigl\|X(\mathbf{e}) \bigr\|_\infty +M_3\bigl\|X(
\mathbf{e})\bigr \|^2_\infty\\
&&{}+M_4\bigl\|X(\mathbf{c})
\bigr\|^2_\infty.
\end{eqnarray*}
Now we prove \eqref{claim}, which holds if the following four
statements hold:
%
\begin{eqnarray}
\mathbb{P} \Bigl(\max_{1\leq|\mathbf{e}-\mathbf{c}| \leq\delta
_n n} M_1\bigl\| X(
\mathbf{e})- X(\mathbf{c})\bigr\|_{\infty}-C_1 \bigl\|V(\mathbf{e})-
\mathbb {I}\bigr\|_1 /16 \leq 0 \Bigr) &\rightarrow&1, \label{s1}
\\
\mathbb{P} \Bigl(\max_{1\leq|\mathbf{e}-\mathbf{c}| \leq\delta
_n n} M_2\bigl\| X(\mathbf{e})
\bigr\|_\infty-C_1 /16 \leq0 \Bigr) &\rightarrow&1, \label{s2}
\\
\mathbb{P} \Bigl(\max_{1\leq|\mathbf{e}-\mathbf{c}| \leq\delta
_n n} M_3\bigl\| X(\mathbf{e})
\bigr\|^2_\infty-C_1\bigl\|V(\mathbf{e})-\mathbb{I}
\bigr\|_1 /16 \leq0 \Bigr) &\rightarrow&1, \label{s3}
\\
\mathbb{P} \Bigl(\max_{1\leq|\mathbf{e}-\mathbf{c}| \leq\delta
_n n} M_4\bigl\| X(\mathbf{c})
\bigr\|^2_\infty- C_1\bigl\|V(\mathbf{e})-\mathbb{I}
\bigr\|_1 /16 \leq0 \Bigr)& \rightarrow&1. \label{s4}
\end{eqnarray}

The proof of \eqref{s1} is similar to the proof of (A.15) in \cite{Zhaoetal2012}.
Note that $\frac{1}{n}|\mathbf{e}-\mathbf{c}| \leq\frac{1}{2} \|
V(\mathbf{e})-\mathbb
{I}\|_1$.
So for each $m \geq1$,
\begin{eqnarray*}
&& \mathbb{P}\Bigl(\max_{|\mathbf{e}-\mathbf{c}|=m} M_1\bigl\|X(\mathbf{e})-
X(\mathbf{c})\bigr\| _{\infty}- C_1\bigl \|V(\mathbf{e})-\mathbb{I}
\bigr\|_1 /16>0 \Bigr)
\\
&&\qquad\leq \mathbb{P}\biggl(\max_{|\mathbf{e}-\mathbf{c}|\leq m} \bigl\|X(\mathbf {e})-X(
\mathbf{c})\bigr\| _{\infty}>\frac{C_1m}{8M_1n}\biggr)=I_1.
\end{eqnarray*}
Let $\alpha=C_1/8 M_1$ if $\alpha\geq6C$, by (A.2) in \cite{Zhaoetal2012},
\begin{eqnarray*}
I_1 & \leq&2 K^{m+2} n^m \exp \biggl( -\alpha
\frac{3m}{8n} \mu_n \biggr)
\\
& =& 2 K^2 \bigl[K \exp\bigl(\log n - \alpha\mu_n/(8/3 n)
\bigr) \bigr]^m.
\end{eqnarray*}
If $\alpha<6C$, by (A.3) in \cite{Zhaoetal2012},
\begin{eqnarray*}
I_1 & \leq&2 K^{m+2} n^m \exp \biggl( -
\alpha^2 \frac{m}{16C n} \mu_n \biggr)
\\
& =& 2 K^2 \bigl[K \exp\bigl(\log n - \alpha^2
\mu_n /(16C n)\bigr) \bigr]^m.
\end{eqnarray*}

In both cases, since $\lambda_n/\log n \rightarrow\infty$ ($\lambda
_n=n \rho_n$),
\begin{eqnarray*}
&&\mathbb{P} \Bigl(\max_{1\leq|\mathbf{e}-\mathbf{c}| \leq\delta
_n n} M_1\bigl\| X(\mathbf{e})-
X(\mathbf{c})\bigr\|_{\infty}-C_1 \bigl\|V(\mathbf{e})-\mathbb {I}
\bigr\|_1 /16> 0 \Bigr)
\\
&&\qquad\leq \sum_{m=1}^{\infty} \mathbb{P} \Bigl(\max
_{|\mathbf
{e}-\mathbf{c}|=m} M_1\bigl\| X(\mathbf{e})- X(\mathbf{c})
\bigr\|_{\infty}-C_1 \bigl\|V(\mathbf{e})-\mathbb {I}\bigr\|_1
/16>0 \Bigr) \rightarrow0,
\end{eqnarray*}
as $n\rightarrow\infty$, which completes the proof of \eqref{s1}.

Equation~\eqref{s2} simply follows (A.1) in \cite{Zhaoetal2012}.

We next prove \eqref{s3}. For each $1\leq m \leq\delta_n n $,
\begin{eqnarray*}
&&\mathbb{P} \Bigl(\max_{|\mathbf{e}-\mathbf{c}|=m} M_3 \bigl\|X(\mathbf {e})
\bigr\|^2_\infty- C_1\bigl\|V(\mathbf{e})-\mathbb{I}
\bigr\|_1 /16 > 0 \Bigr)\\
&&\qquad \leq\mathbb{P} \biggl(\max_{|\mathbf{e}-\mathbf{c}|\leq m}
\bigl\| X(\mathbf{e}) \bigr\| ^2_\infty> \frac{C_1m}{8M_3n}
\biggr)=I_2.
\end{eqnarray*}
Let $\varepsilon=\sqrt{\frac{C_1m}{8M_3n}}$, $\alpha=C_1/64M_3$. Then
from Lemma~\ref{lemma-bernstein},
\begin{eqnarray*}
I_2 & \leq&2K^{m+2} n^m \exp \biggl( -
\frac{3\mu_n \varepsilon
^2}{4(\varepsilon+3)} \biggr)
\\
& \leq&2K^{m+2} n^m \exp \biggl( -\frac{\mu_n \varepsilon^2}{8} \biggr)
\\
& = &2K^{m+2} n^m \exp \biggl(-\alpha\frac{\mu_n}{n} m
\biggr)
\\
& = &2 K^2 \biggl[K \exp \biggl(\log n - \alpha\frac{\mu_n}{n}
\biggr) \biggr]^m. 
\end{eqnarray*}

Since $\lambda_n/\log n \rightarrow\infty$,
\begin{eqnarray*}
&&\mathbb{P} \Bigl(\max_{1\leq|\mathbf{e}-\mathbf{c}| \leq\delta
_n n}M_3 \bigl\| X(\mathbf{e})
\bigr\|^2_\infty- C_1 \bigl\|V(\mathbf{e})-\mathbb{I}
\bigr\|_1 /16> 0 \Bigr)
\\
&&\qquad\leq \sum_{m=1}^{\infty} \mathbb{P} \Bigl(\max
_{|\mathbf
{e}-\mathbf{c}|=m} M_3\bigl \|X(\mathbf{e}) \bigr\|^2_\infty-C_1
\bigl\|V(\mathbf{e})-\mathbb{I}\bigr\|_1 /16> 0 \Bigr) \rightarrow0,
\end{eqnarray*}
as $n\rightarrow\infty$, which completes the proof of \eqref{s3}.

We now complete the proof by showing \eqref{s4}.
For each $1\leq m \leq\delta_n n $,
\begin{eqnarray*}
&&\mathbb{P} \Bigl(\max_{|\mathbf{e}-\mathbf{c}|=m} M_4\bigl \|X(\mathbf {c})
\bigr\|^2_\infty - C_1\bigl\| V(\mathbf{e})-\mathbb{I}
\bigr\|_1 /16 > 0 \Bigr) \\
&&\qquad= \mathbb{P} \biggl( \bigl\| X(\mathbf{c}) \bigr\|
^2_\infty> \frac{C_1m}{8M_4 n} \biggr)=I_3.
\end{eqnarray*}
Let $\varepsilon=\sqrt{\frac{C_1m}{8M_4n}}$, $\alpha=C_1/64M_4$. Then from
Bernstein's inequality,
%
\begin{equation}
I_3  \leq2K^2 \exp \biggl( -\frac{3\mu_n \varepsilon^2}{4(\varepsilon+3)} \biggr)
\leq2K^2 \exp \biggl(-\alpha\frac{\mu_n}{n} m \biggr).
\end{equation}
Therefore,
\begin{eqnarray*}
&&\mathbb{P} \Bigl(\max_{1\leq|\mathbf{e}-\mathbf{c}| \leq\delta
_n n}M_4 \bigl\| X(\mathbf{c})
\bigr\|^2_\infty-C_1 \bigl\|V(\mathbf{e})-\mathbb{I}
\bigr\|_1 /16> 0 \Bigr)
\\
& &\qquad\leq\sum_{m=1}^{\infty} \mathbb{P}\bigl(
M_4 \bigl\|X(\mathbf{e})\bigr \| ^2_\infty-
C_1 \bigl\|V(\mathbf{e})-\mathbb{I}\bigr\|_1 /16> 0 \bigr)
\rightarrow0 \qquad\mbox{as } n\rightarrow\infty.
\end{eqnarray*}


\section*{Acknowledgements}
We are very grateful to Emma Jingfei Zhang, a former Ph.D. student at
University of Illinois at Urbana-Champaign now at University of Miami,
who discovered the error and
persisted in tracking down its root cause.

%



\printaddresses